\documentclass[12pt,twoside]{article}
\usepackage{amsmath,amssymb,mathrsfs,graphicx,subfigure}
\usepackage[small,sc]{caption} 
\usepackage[nodvipsnames]{color}
\usepackage{cancel}
\usepackage[normalem]{ulem}

\newtheorem{theorem}{Theorem}[section]
\newtheorem{lemma}[theorem]{Lemma}
\newtheorem{proposition}[theorem]{Proposition}

\newtheorem{rmrk}[theorem]{Remark}

\DeclareMathAlphabet{\mathbfit}{OML}{cmm}{b}{it}

\makeatletter
\@addtoreset{equation}{section}
\makeatother

\newenvironment{remark}
{\begin{rmrk} \em}
{\end{rmrk}}

\newcommand{\sy} {system}

\newcommand{\R} {\mathbb{R}}

\newcommand{\Q} {\mathbb{Q}}
\newcommand{\Z} {\mathbb{Z}}
\newcommand{\N} {\mathbb{N}}
\newcommand{\HHH} {\mathbb{H}}
\newcommand{\qed} {\hfill {\small Q.E.D.} \par\medskip}

\newcommand{\ds} {\displaystyle}
\newcommand{\proof} {\noindent \textsc{Proof.} }
\newcommand{\proofof}[1] {\noindent \textsc{Proof of {#1}.} }
\newcommand{\article}[3] {\textsc{{#1}}, {\itshape {#2}}, {{#3}}.}
\newcommand{\book}[3] {\textsc{{#1}}, {\itshape {#2}}, {{#3}}.}
\newcommand{\vol} {\textbf}
\newcommand{\eps} {\varepsilon}

\newcommand{\rset}[2] {\left\{ #1 \: \left| \: #2 \right. \! \right\} }

\newcommand{\into} {\longrightarrow}

\newcommand{\newfig}[4] {
\bigskip\medskip
\begin{figure}[htbp]
  \centering
  \includegraphics[width=#3]{#2}
  \begin{minipage}[t]{0.80\linewidth} 
    \caption{#4}
    \protect\label{#1}
  \end{minipage}
\end{figure}
\smallskip
}

\newcommand{\btable} {\Omega}
\newcommand{\map} {f}
\newcommand{\torus} {\mathbb{S}^1}
\newcommand{\alpham} {\alpha_M}  
\newcommand{\ubtable} {\btable_4}  
\newcommand{\dilsurf} {\btable'}  
\newcommand{\perset} {\mathcal{P}}
\newcommand{\id} {\mathrm{Id}}

\begin{document}

\title{\textbf{Internal-wave billiards in \\ trapezoids and similar tables}}

\author{
\scshape
Marco Lenci\thanks{
Dipartimento di Matematica, Universit\`a di Bologna,
Piazza di Porta San Donato 5, 40126 Bologna, Italy. 
E-mail: \texttt{marco.lenci@unibo.it}.}
\thanks{
Istituto Nazionale di Fisica Nucleare,
Sezione di Bologna, Viale Berti Pichat 6/2,
40127 Bologna, Italy.}
,
Claudio Bonanno\thanks{
Dipartimento di Matematica, Universit\`a di Pisa, Largo Bruno 
Pontecorvo 5, 56127 Pisa, Italy. E-mail: 
\texttt{claudio.bonanno@unipi.it}.}
,
Giampaolo Cristadoro\thanks{
Dipartimento di Matematica e Applicazioni, Universit\`a di Milano - Bicocca,
Via Roberto Cozzi 55, 20125 Milano, Italy. E-mail: 
\texttt{giampaolo.cristadoro@unimib.it}.}
}

\date{Final preprint for \emph{Nonlinearity} \\[6pt]
  December 2022}

\maketitle

\begin{abstract}
  We call \emph{internal-wave billiard} the dynamical system of a point
  particle that moves freely inside a planar domain (the \emph{table}) and is 
  reflected by its boundary according to this nonstandard rule: the 
  angles that the incident and reflected velocities form with a fixed direction 
  (representing gravity) are the same. These systems are point particle 
  approximations for the motion of internal gravity waves in closed containers, 
  hence the name. For a class of tables similar to rectangular trapezoids, 
  but with the slanted leg replaced by a general curve with downward concavity,
  we prove that the dynamics has only three asymptotic regimes: (1) there exist 
  a global attractor and a global repellor, which are periodic and might coincide; 
  (2) there exists a beam of periodic trajectories, whose boundary (if any) 
  comprises an attractor and a repellor for all the other trajectories; (3) all 
  trajectories are dense (that is, the system is minimal). Furthermore, in the 
  prominent case where the table is an actual trapezoid, we study the sets in 
  parameter space relative to the three regimes. We prove in particular that the 
  set for (1) has positive measure (giving a rigorous proof of the existence of 
  Arnol'd tongues for internal-wave billiards), whereas the sets for (2) and (3) 
  are non-empty but have measure zero.
\end{abstract}

\section{Introduction}
\label{sec-intro}

When an incompressible fluid is stably stratified by a linear increase 
of density in the direction of gravity, a periodic perturbation can generate 
gravity waves whose direction of propagation is dictated only by the 
frequency of the forcing \cite{Turner73}. These \emph{internal waves}
appear ubiquitously in oceans and in the atmosphere, where they play a 
key role in mixing processes and energy dissipation \cite{s}. Sometimes
the stratification in density is due to a centrifugal force or other
inertial effects. In such cases the resulting internal waves are also called
\emph{inertial waves} \cite{Sibgatullin2019}.

When internal waves encounter a rigid body, they are reflected according 
to this nonstandard rule \cite{Maas2005, Pillet2019}:
\begin{itemize}
\item[(a)] the angles that the incident and reflected velocities make with 
  the vertical direction are equal or supplementary (in other words, the 
  two vectors belong to the boundary of the same vertical cone);

\item[(b)] the orthogonal projections of the incident and reflected 
  velocities on the normal to the body at the reflection point are
  opposite to each other;
  
\item[(c)] the orthogonal projections of the incident and reflected 
  velocities on the direction orthogonal to both the normal and the
  vertical direction are the same;

\item[(d)] if, given an incident velocity, the conditions (a)-(c) do not 
  determine a unique reflected velocity (which is the general case), 
  out of the two vectors which satisfy (a)-(c), the reflected velocity is 
  the one whose projection on the plane generated by the normal and 
  the vertical is \emph{not} the opposite of the projection of the incident 
  velocity on the same plane.
\end{itemize}
Observe that the standard (Fresnel) reflection rule is the same as the
above except that (a) is replaced by the condition that the angles of 
incidence and of reflection \emph{relative to the normal} are the same,
and (d) becomes redundant.
It can be checked that with the above rule, which we refer to as the 
\emph{internal-wave reflection rule}, the moduli of the incident and 
reflected velocities are in general not the same \cite[Sect.~2.2]{Pillet2019} 
--- with the Fresnel rule they are.

It is also easy to see that the internal-wave reflection rule causes a beam 
of parallel rays to expand and/or contract its section, in general, at every 
reflection. In their seminal paper \cite{Maas1995}, Maas and Lam 
predicted that, in a closed container, the contraction effect
prevails in certain parts of phase
space so as to generate a periodic attractor for the waves. This 
\emph{internal wave attractor}, as it was dubbed, is determined solely 
by the geometry of the container and the frequency of the periodic 
perturbation. This prediction was soon confirmed experimentally 
\cite{Maas1997} and later other laboratories \cite{Hazewinkel2008, 
Hazewinkel2010, Hazewinkel2011, bejpd} and numerical simulations 
\cite{Grisouard2008, bssed} verified and extended the findings of 
\cite{Maas1997}. (Here and in the rest of this paper we do not pretend
to cite the entire literature on the subject of internal waves, which is 
massive: we only reference the publications that are closest to our 
work.)

From a theoretical point of view, the first step towards understanding 
the structure of internal waves in closed containers consists in studying 
its \emph{ray dynamics} \cite{Maas1995, Maas1997, Manders2003b,
Maas2005, Lam2008, Hazewinkel2008, Hazewinkel2010, Pillet2019,cs},
which is defined by the free motion of a point particle inside the container, 
subject to the internal-wave reflection law at its boundary. This dynamics 
has been studied mostly in two-dimensional domains, which one thinks 
of as vertical, relative to gravity.
In this case the velocity of the particle can only assume four directions.
If we indicate a direction by means of the angle $\vartheta$ it forms with 
the direction of gravity, the four possibilities are: $\theta, \pi-\theta, 
\pi+\theta, -\theta$, where $\theta$ is an angular parameter that is related 
to the external forcing in the original wave \sy. (In this paper $\theta$ will 
play the equivalent role of a certain initial direction of the motion, see below.) 

We denote the domain by $\btable$ and assume that its boundary
$\partial \btable$ is piecewise smooth. When the particle hits a point $P \in
\partial \btable$, its velocity changes instantaneously and takes the unique 
direction, among the 4 possibilities, that points towards the interior of 
$\btable$ and is not the opposite of the incoming direction. (This rule is 
ambiguous when the slope of $\partial \btable$ at $P$ is either undefined 
or in one of the four special directions. We discard these trajectories, which 
amount to a Lebesgue zero-measure set in four-dimensional phase space.) 
As for the modulus of the outgoing velocity, this is a function of the incoming 
velocity and the slope of $\partial \btable$ at $P$ \cite[Sect.~2.2]{Pillet2019}:
we do not recall it here because it will soon become irrelevant. Thus, there 
are two types of reflections: \emph{vertical} reflections, where the direction 
of the velocity turns from a given $\vartheta$ to $-\vartheta$, and 
\emph{horizontal} reflections, where the direction turns from $\vartheta$ 
to $\pi-\vartheta$. 

We call \emph{internal-wave billiard} any system like the above, 
with the difference that the speed of the particle, i.e, the modulus 
of the velocity, \emph{is constantly equal to 1}. This reparametrization of 
time naturally changes the dynamical system, so the ray dynamics in 
$\btable$ is \emph{not} the internal-wave billiard in $\btable$, but if one is 
only interested in combinatorial/topological properties (say, 
periodicity or density of the trajectories, attractors, etc.), the two 
systems are equivalent. 

Internal-wave billiards are special cases of a larger class of billiard-like
systems sometimes called \emph{chess billiards}; see \cite{hm, nt} and
references therein. In the rest of the paper we only deal with 
internal-wave billiards, referring to $\btable$ as the \emph{(billiard) table}. 

Different tables have been investigated 
\cite{Maas1995, Manders2003b, bfm}, showing the ubiquity of 
internal wave attractors, but the rectangular trapezoid, which was the 
shape chosen in the first actual experiment with internal waves 
\cite{Maas1997}, emerged as a prototypical example. When 
$\btable$ is a rectangular trapezoid with horizontal bases and a 
vertical leg, reflections are always of horizontal type at the bases and of 
vertical type at the vertical leg. As for the slanted leg, let us look at 
Fig.~\ref{fig1}\emph{(b)} to define the angles $\alpha$ and $\theta$, 
the latter denoting the initial direction of the particle, conventionally taken 
at the vertical leg (we may also assume $\theta \in [0, \pi/2]$, which
is no loss of generality, as will be made clear later): reflections at
the slanted leg are horizontal for $0 < \theta < \alpha$ and vertical for 
$\alpha < \theta < \pi/2$. We do not consider the cases $\theta \in 
\{0, \alpha, \pi/2\}$ because they are formally ill-defined and/or trivial. 
The case $0 < \theta < \alpha$ is also simple: all trajectories converge 
to the rightmost corner of the trapezoid. It is for $\alpha < \theta < \pi/2$ 
that the dynamics gets interesting: previous work suggests that only 
three scenarios are possible: 
(1) there exist a unique global attractor and a unique global repellor 
(a.k.a.\ backward attractor), which might possibly coincide: they are periodic 
with the same number of bounces off the boundary of $\btable$;
(2) all trajectories are periodic with the same number of bounces off the 
boundary;
(3) all trajectories are dense.
While the first and second scenarios have been routinely observed in 
the physical literature, little is known about the occurrence of the third one.

\newfig{fig1}{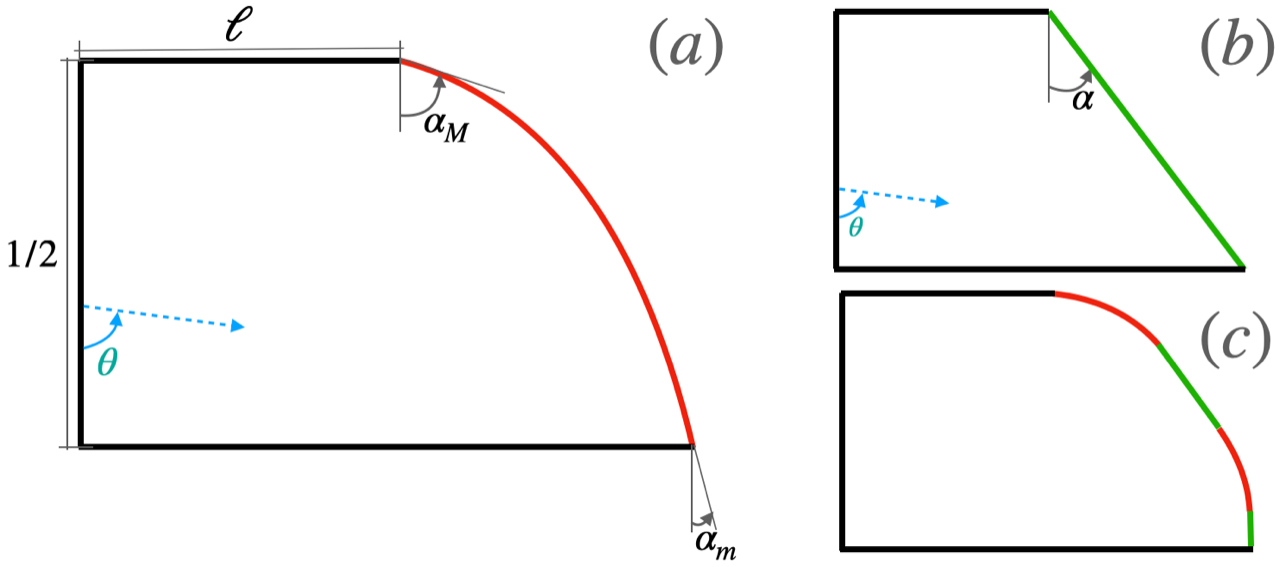}{11cm}{Examples of billiard tables $\btable$.} 

In this paper we give rigorous results for the internal-wave billiards
in a class of tables which includes the rectangular trapezoid. 
Specifically, we allow the slanted leg of the trapezoid to be replaced 
by a piecewise smooth curve with downward concavity; see 
Fig.~\ref{fig1} for examples of tables in this class and Section 
\ref{sec-1ddyn} for a precise definition of it. Our main result gives a 
comprehensive description of the asymptotic 
behavior of all trajectories of the system. In 
particular, it is proven that there are only three possibilities, with 
subcases: 
\begin{itemize}
\item[(1)] There exist a unique global attractor and a unique global repellor, 
which are periodic with the same number of bounces off the boundary
of the table, and might possibly coincide.

\item[(2)] There exists a beam of periodic trajectories. 
  \begin{itemize}
  \item[(a)] If the beam is non-degenerate, its boundary consists of two 
  periodic trajectories, which are, respectively, the attractor and the 
  repellor of all trajectories outside of the beam; otherwise
  
  \item[(b)] the beam may reduce to a single periodic trajectory, which  
  acts as both a global attractor and a global repellor; or

  \item[(c)] the beam may comprise all trajectories, resulting in the case 
  where all trajectories are periodic with the same number of bounces 
  off the boundary.
  \end{itemize}

\item[(3)] All trajectories are dense.
\end{itemize}

Furthermore, if the slanted side of $\btable$ does not contain any segment, 
case (2) coincides with subcase (b); if the slanted side of $\btable$ is a 
segment (i.e., $\btable$ is a rectangular trapezoid), case (2)  coincides with 
subcase (c). All these statements are direct consequences of Theorem 
\ref{main-thm}, which describes the dynamics of a suitable Poincar\'e map 
for the billiard. Section \ref{sec-1ddyn} is devoted to introducing, stating 
and proving Theorem \ref{main-thm}.

In Section \ref{sect-trapez} we study the trapezoidal case in depth, 
finding in particular that, in parameter space, case (1) has full measure, and 
case (2) (which is the same as subcase (c) here) and case (3) are non-empty 
with zero measure. The statement about case (1) is important, we believe, 
because it is the first mathematical proof --- to our knowledge --- of the 
existence of Arnol'd tongues (definition in Section \ref{sect-trapez}) for 
internal-wave billiards. The statement about case (3) is also interesting, 
because it explains why this regime was hardly observed in experiments and 
numerical simulations.

\paragraph{Acknowledgments.} We thank Thierry Dauxois for useful
discussions in the early stage of this work and Selim Ghazouani for helping 
us correct a mistake in an earlier version of this manuscript. The present 
research was partially supported by the PRIN Grant 2017S35EHN, MUR, 
Italy. It is also part of the authors' activity within the UMI Group 
\emph{DinAmicI} and the Gruppo Nazionale di Fisica Matematica, INdAM.

\section{Reduction to one-dimensional dynamics}
\label{sec-1ddyn}

We study the internal-wave billiard in a table $\btable$ that is a 
generalization of a rectangular trapezoid, like the examples shown in 
Fig.~\ref{fig1}. In suitable units, $\btable$ has height 1/2 and is 
completely specified by the position and shape of its rightmost 
boundary $\partial_R \btable$, which we assume to be the graph of a 
piecewise $C^1$, strictly decreasing, concave function, possibly with
the addition of a vertical segment attached to its lower end. 
Fig.~\ref{fig1}\emph{(a)} defines three 
important parameters for the dynamics, the angles $\alpha_m$, 
$\alpham$ and $\theta$. In particular, $\theta \in (0,\pi/2) \cup 
(\pi/2, \pi)$ is the particle's initial direction, having 
assumed (without loss of generality) an initial position on 
$\partial_L \btable$, the vertical side of $\btable$. We do not consider 
horizontal or vertical initial directions for the reasons explained 
earlier.

If $\theta \in (0,\alpha_m) \cup (\pi - \alpha_m, \pi)$, it is easily verified
that all reflections at $\partial_R \btable$ are horizontal (as per the
definition given in the introduction), implying that every trajectory 
converges to the lower endpoint of $\partial_R \btable$. If
$\theta \in [\alpha_m, \alpham] \cup [\pi - \alpham, \pi - \alpha_m]$,
$\partial_R \btable$ is generally split into two parts, which give rise to 
horizontal and vertical reflections, respectively. This is the more 
complicated case, which we do not study at the present time, as we 
intend to keep our mathematical machinery to a minimum. Moreover,
when $\btable$ is a trapezoid, which a case of interest here, the 
dynamics for $\theta = \alpha$ is ill-defined.

In this paper we restrict to the case $\theta \in (\alpham, \pi/2) \cup 
(\pi/2, \pi-\alpham)$ to ensure that all reflections at $\partial_R \btable$
are vertical. Since the particle's velocity can only assume 4 values, 
it is convenient to represent the dynamics on $\btable$ as a 
\emph{linear flow} on $\ubtable$, the 4-fold copy of $\btable$ 
represented in Fig.~\ref{fig2} as a subset of $\R^2$. The linear flow 
on $\ubtable$ is defined as a fixed-velocity motion on it, with the 
provision that a trajectory hitting a boundary point continues on 
the opposite boundary point with the same velocity. Two boundary 
points of $\ubtable$ are said to be opposite to each other if they 
belong to the upper/lower boundaries of $\ubtable$ and have the 
same abscissa, or to the left/right boundaries and have the same 
ordinate. Fig.~\ref{fig2} also displays a trajectory of the linear flow, 
together with its projection on $\btable$. The defining parameter of 
the linear flow is the angle $\theta$ introduced earlier, 
which we also call the \emph{direction} of the flow. We convene that
the speed of the flow is 1. It is apparent that the natural projection 
$\ubtable \into \btable$ maps trajectories of the flow into trajectories 
of the billiard. The procedure whereby one passes from the 
internal-wave billiard on $\btable$ to the linear flow on $\ubtable$ is 
also called \emph{billiard unfolding}. (This technique has been applied
very fruitfully to polygonal (ordinary) billiards; see, e.g., the reference 
list of \cite{ddl}.)

\newfig{fig2}{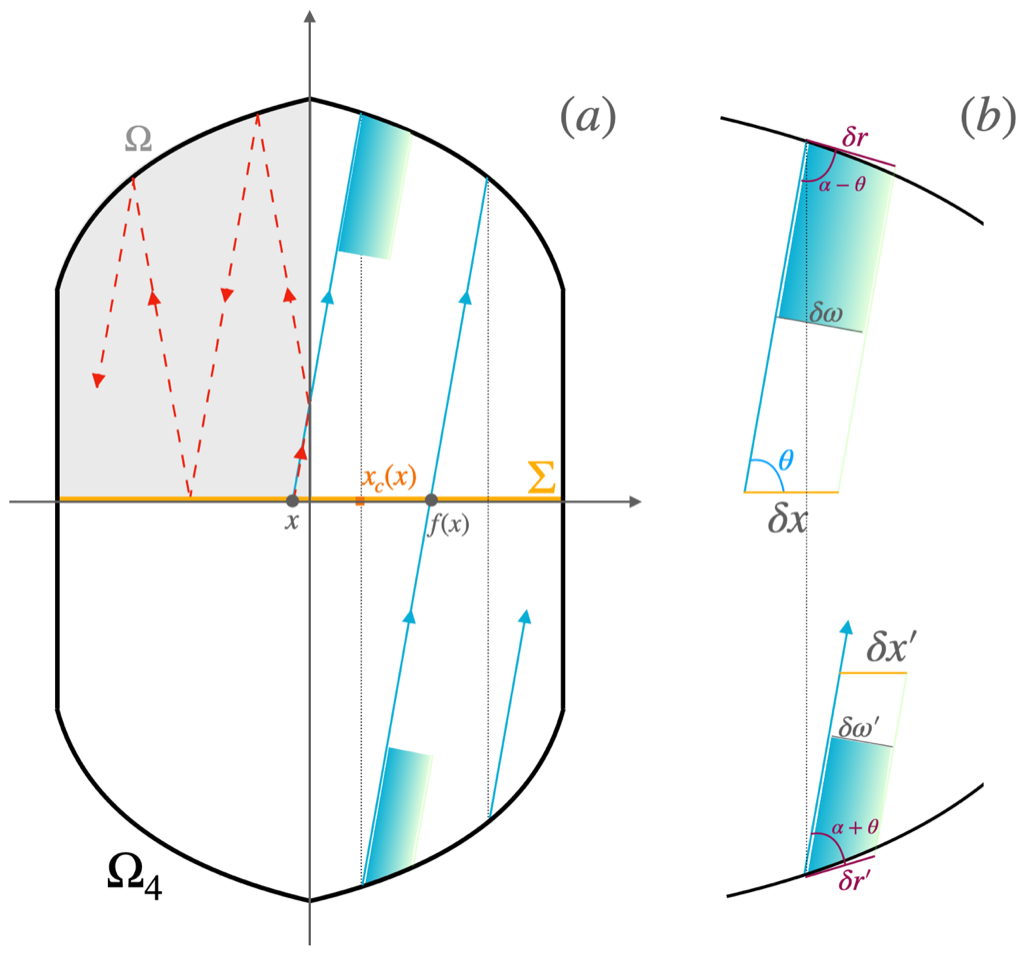}{10cm}{A trajectory of the linear flow in 
$\ubtable$, together with the corresponding orbit for the Poincar\'e 
map $\map: \Sigma \into \Sigma$. The definition of $x_c(x)$, cf.\ 
Lemma \ref{lem-xc}, is also illustrated.} 

Our assumptions on the billiard table $\btable$, given at the 
beginning of this section, can be restated as assumptions on 
$\ubtable$ as follows: the upper boundary of $\ubtable$ is the 
graph of a function $b: [-1/2, 1/2] \into \R^+$, which is even, 
piecewise $C^1$, concave, and (not necessarily strictly) decreasing 
on $[0,1/2]$. Set $\alpha(x) := \arctan(b'(x))$. This expression fails to 
be defined in at most countably many points of $(-1/2, 1/2)$. The 
assumptions show that $\alpha$ is a (not necessarily strictly) 
decreasing odd function with $\alpha(-1/2) = \alpham$ and 
$\alpha(1/2) = -\alpham$. In particular, when $b'(x)$ is defined, 
$x \le 0$ implies that $\alpha(x) \ge 0$ and $x \ge 0$ implies that 
$\alpha(x) \le 0$.

Identifying opposite boundary points of $\ubtable$, which we 
henceforth do, we regard $\ubtable$ as a surface homeomorphic to 
a torus. The metric of $\R^2$ induces on this surface a metric 
which is defined and flat everywhere except for a closed curve (which, 
in topological terms, is a simple, non-contractible loop). The linear 
flow is the geodesic flow for this metric, subject to a non-isometric 
identification rule between the two sides of the closed curve.

A convenient way to study the elementary properties of a flow is by 
means of a suitable Poincar\'e map. In our case, a good choice is to 
take the first-return map to the horizontal segment $\Sigma$ shown 
in Fig.~\ref{fig2}\emph{(a)}. We denote it $\map: \Sigma \into \Sigma$. 
As per our 
boundary identifications, $\Sigma$ has the topology of a circle. In any 
case, with the units that we have chosen, its length is 1. Observe that 
points on this Poincar\'e section correspond to the particle being on 
$\partial_L \btable$ with velocity directed as $\theta$ or $\pi-\theta$. 
In other words, $\Sigma$ is in 2-to-1 correspondence with the set 
of initial positions we chose for our dynamics. This fact and the 
symmetry of $\ubtable$ show that it is no loss of generality to restrict 
the directions of the linear flow to $\theta \in (\alpham, \pi/2)$.  
Identifying $\Sigma$ with $\torus := \R/\Z$ in the natural way, we also 
regard $\map$ as a homeomorphism of $\torus$. 

\begin{lemma} \label{lem-xc}
  Let us denote by $x_c(x)$ the abscissa of the first collision point of 
  the trajectory starting from $x \in \Sigma$ with the upper boundary of 
  $\ubtable$ (see Fig.~\ref{fig2}(a)). For all $x \in \Sigma$ such that 
  $\alpha(x_c(x))$ is defined, we have
  \begin{displaymath}
    \map'(x) = \frac{\sin(\theta + \alpha(x_c(x)))} 
    {\sin(\theta-\alpha(x_c(x)))} .
  \end{displaymath}
\end{lemma}

\proof This proof will be illustrated by Fig.~\ref{fig2}\emph{(b)}. Given
$x \in \Sigma$ such that $\alpha(x_c(x))$ is well defined in a neighborhood 
of $x$, let us take $\delta x$
so small that such neighborhood contains $[x, x + \delta x]$ 
(having chosen $\delta x > 0$ without loss of generality). The beam of
(parallel) trajectories stemming from $[x, x + \delta x]$ has 
width $\delta w = \delta x \sin \theta$ and thus projects on the 
upper boundary of $\ubtable$ an arc of length 
\begin{equation}
  \delta r = \frac{\delta w} {\sin(\theta - \alpha)} + o(\delta x) =
  \delta x \, \frac{\sin \theta} {\sin(\theta - \alpha)} + o(\delta x), \quad
  \mbox{ as } \delta x \to 0.
\end{equation}
Here, with a slight abuse of notation, we have denoted 
$\alpha := \alpha(x_c(x))$. The beam then continues from
the opposite arc, on the lower boundary of $\ubtable$. The length of
this arc is clearly $\delta r' = \delta r$ but, reversing the previous
reasoning, its width is $\delta w' = \delta r' \sin(\theta + \alpha) + 
o(\delta x)$. The first intersection of this beam with $\Sigma$ is the 
segment $[\map(x), \map(x + \delta x)]$, whose length is
\begin{equation}
  \map(x + \delta x) - \map(x) = \frac{\delta w'} {\sin \theta}  = \delta x \, 
  \frac{\sin (\theta + \alpha)} {\sin(\theta - \alpha)} + o(\delta x), \quad
  \mbox{ as } \delta x \to 0.
\end{equation}
Dividing the above by $\delta x$ and taking the limit proves the lemma.
\qed

Another property of $\map$ is that, for $x, x_1 \in \Sigma$,
\begin{equation}
  x_1 = \map(x) \quad \Longleftrightarrow \quad -x = \map(-x_1).
\end{equation}
In other words, the graph of $\map$, as represented in $\R^2$, is
symmetric around the bisectrix of the second and fourth quadrants.
This is easily seen by drawing a first-return segment of a trajectory of 
the linear flow on $\ubtable$ and rotating $\ubtable$ by 180 degrees, 
equivalently, by exploiting the fact that the billiard dynamics commutes
with time-inversion.

In order to study the topological/combinatorial properties of the 
internal-wave billiard, it suffices to study the corresponding properties 
of $\map$ on $\torus$. There exists a well-developed theory of circle 
homeomorphisms, dating back from Poincar\'e, who introduced the 
notion of rotation number, which we briefly recall. Let $F: \R \into \R$ 
be a \emph{lift} of $\map$, i.e., a homeomorphism of $\R$ which is well 
defined on the equivalence classes of $\R/\Z$ (this turns out to be the 
same as the property $F(x + k) = F(x) + k$, for all $x\in\R$ and 
$k\in\Z$) and acts like $\map$ there. Then the quantity
\begin{equation} \label{rot-no}
  \rho(\map) := \lim_{n \to \infty} \frac{F^n(x) - x} n
\end{equation}
does not depend on $x$. If taken mod 1, as is customary, it is also
independent of the choice of $F$. One calls $\rho(\map)$ the 
\emph{rotation number} of $\map$. In this paper, unless otherwise 
stated, we represent it as a number in $[0,1)$. Standard references 
in this area are, e.g., \cite[Chap.~1]{dv} and \cite[\S 11]{kh}. Using only 
basic results, the next theorem describes the asymptotic behavior of 
\emph{all} the orbits of $\map$. 

\paragraph{Notation.} In what follows we will need to use 
\emph{cyclic indices}, namely, elements of $\Z_q := \Z/q\Z$, for some 
$q \in \Z^+$. We can think of them as elements of $\{ 0, 1, \ldots, q-1 \}$ 
with the understanding that the sum of any such element with an integer 
is intended mod $q$.
\medskip

\begin{theorem} \label{main-thm}
  Let $\rho(\map) \in [0,1)$ denote the rotation number of $\map: \torus
  \into \torus$.
  \begin{itemize}
  \item[(i)] If $\rho(\map) \not\in \Q$, $\map$ is topologically equivalent
    (i.e., conjugated via homeomorphism) to the rotation $x \mapsto x + 
    \rho(\map) \ (\mathrm{mod}\ 1)$ on $\torus$. In particular, $\map$ 
    is minimal, i.e., all orbits are dense. 
    
  \item[(ii)] If $\rho(\map) \in \Q$, we write $\rho(\map) = p/q$, where, if 
    $\rho(\map)>0$, $p$ and $q$ are coprime positive integers, or else 
    $p=0$ and $q=1$. There exist two sets
    \begin{displaymath}
      \perset^\pm := \bigcup_{i=0}^{q-1} \perset_i^\pm, 
    \end{displaymath}
    where the $\perset_i^\pm$ are closed intervals. All points in $\perset^\pm$ 
    are $q$-periodic and such that, for all $i \in \Z_q$ and $n \in \N$,
    \begin{displaymath}
      \map^n \!\left( \perset_i^\pm \right) = \perset_{i+pn}^\pm .
    \end{displaymath}
    All points in $\torus \setminus ( \perset^+ \cup \perset^-)$ are non-periodic. 
    
    Only two possibilities are given:
    \begin{itemize}
    \item[(a)] Each $\perset_i^\pm$ reduces to a point, which we denote 
      $x_i^\pm$. These points are ordered as follows, according to the 
      orientation of $\torus$:
      \begin{displaymath}
        x_0^- \le x_0^+ \le x_1^- \le x_1^+ \le \cdots \le x_{q-1}^- \le x_{q-1}^+ 
        \le x_0^- .
      \end{displaymath}
      So $\mathcal{O}^\pm := \perset^\pm = \{ x_0^\pm, x_1^\pm, \ldots, 
      x_{q-1}^\pm \}$ are two (possibly coinciding) periodic orbits such that 
      $\map^n (x_i^\pm) = x_{i+pn}^\pm$, for all $i \in \Z_q$ and $n \in \N$. As for 
      the other points:
      \begin{align*}
        \forall x \in \left( x_i^- \,, x_i^+ \right), \qquad & \map^n(x) \in 
        \left( x_{i+pn}^- \,, x_{i+pn}^+ \right) \mbox{ and} \\[3pt]
        & \lim_{n \to +\infty} \map^{nq}(x) = x_i^+ \,, \quad
        \lim_{n \to -\infty} \map^{nq}(x) = x_i^- \,; \\[6pt]
        \forall x \in \left( x_i^+ \,, x_{i+1}^- \right), \qquad & \map^n(x) \in
        \left( x_{i+pn}^+ \,, x_{i+1+pn}^- \right) \mbox{ and} \\[3pt]
        & \lim_{n \to +\infty} \map^{nq}(x) = x_i^+ \,, \quad
        \lim_{n \to -\infty} \map^{nq}(x) = x_{i+1}^- \,.
      \end{align*}
      Hence $\mathcal{O}^+$, respectively $\mathcal{O}^-$, is the global 
      attractor, respectively repellor, of the \sy.
 
    \item[(b)] All $\perset_i^\pm$ are non-degenerate intervals, with
      $\perset_i^+ = \perset_i^-$, for all $i \in \Z_q$. Denoting any such 
      interval $[x_i^L, x_i^R] := \perset_i^+ = \perset_i^-$, we have 
      \begin{displaymath}
        x_0^L < x_0^R \le x_1^L < x_1^R \le \cdots \le x_{q-1}^L < 
        x_{q-1}^R \le x_0^L ,
      \end{displaymath}
      with the property that $x_i^R < x_{i+1}^L$ holds for some $i \in \Z_q$ 
      if, and only if, it holds for all $i \in \Z_q$. 
      \begin{itemize}
      \item[(1)] If $x_i^R < x_{i+1}^L$ for all $i$ then, for all $i \in \Z_q$ and 
        $n \in \N$, $\map^n (\perset_i^\pm) = \perset_{i+pn}^\pm$ and
        \begin{align*}
          \forall x \in \left( x_i^R \,, x_{i+1}^L \right), \qquad & \map^n(x) \in
          \left( x_{i+pn}^R \,, x_{i+1+pn}^L \right) \mbox{ and} \\[3pt]
          & \lim_{n \to +\infty} \map^{nq}(x) = x_i^R \,, \quad
          \lim_{n \to -\infty} \map^{nq}(x) = x_{i+1}^L \,.
        \end{align*}
        Therefore $\mathcal{O}^+ := \{ x_0^R, x_1^R, \ldots, x_{q-1}^R \}$
        and $\mathcal{O}^- := \{ x_0^L, x_1^L, \ldots, x_{q-1}^L \}$ are,
        respectively, the unique (but not global) attractor and repellor of the 
        \sy.
    
     \item[(2)] If $x_i^R = x_{i+1}^L$ for all $i$, all orbits are periodic.
     \end{itemize}
  \end{itemize}
    
  Under the additional assumption that $\partial_R \btable$ does 
  not contain any segment (equivalently, the function $b$ whose graph is 
  the upper boundary of $\ubtable$ is strictly concave), case (b) cannot 
  occur.
  \end{itemize}
\end{theorem}

\begin{remark}
  Since an orbit of $\map$ is dense/periodic/attracting/repelling if, 
  and only if, the corresponding flow trajectory in $\ubtable$ is, all the 
  statements of Theorem \ref{main-thm} are immediately translated to 
  statements about the internal-wave billiard in $\btable$. In particular,
  this proves the claim made in the introduction about the three sole
  possibilities for the asymptotics of the trajectories.
\end{remark}

\proofof{Theorem \ref{main-thm}} Poincar\'e's classical theory of 
circle homeomorphisms \cite[Sect.~I.1]{dv} states that, for any
orientation-preserving $\map: \torus \into \torus$, the following dichotomy 
holds:
\begin{enumerate}
\item If no periodic orbit exists, then $f$ is topologically semiconjugate 
  to an irrational rotation \cite[Thm.~1.1]{dv}, which must necessarily be 
  the rotation by $\rho(\map) \not\in \Q$.
  
\item If a periodic orbit exists and $q$ is its (primitive) period, 
  then necessarily $\rho(\map) = p/q$, where $p$ is either 0 or coprime to $q$. 
  Let $x_0, x_1, \ldots, x_{q-1}$ be an orientation-preserving labeling of the 
  points of the periodic orbit. By definition of rotation number, for all $i \in \Z_q$
  and $n \in \N$, one has
  \begin{equation} \label{main-10}
    \map^n(x_i) = x_{i+pn}
  \end{equation}
  (recall the convention on cyclic indices). Therefore 
  \begin{equation} \label{main-20}
    \map^n \!\left( [x_i, x_{i+1}) \right) = [x_{i+pn}, x_{i+1+pn}).
  \end{equation}
  Necessarily, then, if $x'$ is another periodic point, its combinatorics is the
  same as that of $x$, that is, its period is $q$ and, for any 
  orientation-preserving labeling $x'_0, x'_1, \ldots, x'_{q-1}$ of its orbit, the 
  analogue of (\ref{main-10}) holds. 
\end{enumerate}

To establish \emph{(i)} it suffices to verify that, for our particular $\map$, 
the topological semiconjugacy of case 1 is in fact a topological conjugacy. 
This is exactly the assertion of Denjoy's Theorem for certain circle 
diffeomorphisms \cite[Sect.~I.2]{dv}. As explained in \cite[Rmk on p.~38]{dv},
Denjoy's Theorem also holds for homeomorphisms $\map$ that are 
piecewise differentiable and such that $\log |\map'|$ can be extended to a 
map with bounded variation. Our $\map$ falls in this category, as shown
momentarily. 

Let us call \emph{break point of $\map'$} any $x \in \torus$ such that 
\begin{equation}
  \Delta \map'(x) := \lim_{s \to x^+} \map'(s) - \lim_{s \to x^-} \map'(s)
  \ne 0,
\end{equation}
provided the limits exist. A break point $x$ of $\map'$ is said to be 
\emph{of increase} or \emph{decrease} if $\Delta \map'(x) > 0$ or
$\Delta \map'(x) < 0$, respectively.

\begin{lemma} \label{lem-non-reg}
  The homeomorphism $\map$ defined earlier is such that 
  $\lim_{s \to x^\pm} \map'(s)$ exists at all $x \in \torus$; $\map'$ has 
  exactly one break point of increase, denoted $a_0$, and at most 
  countably many break points of decrease, denoted $\{a_i\}_{i \ge 1}$. 
  This implies that $\map'$ is continuous on 
  $\torus \setminus \{a_i\}_{i\ge 0}$ with positive one-sided limits 
  everywhere. Furthermore, $\map$ is concave on the arc 
  $\torus \setminus \{a_0\}$. Finally, $\log \map'$ can be extended to a 
  map with bounded variation.
\end{lemma}

\proofof{Lemma \ref{lem-non-reg}} By construction of the function
$b: \torus \into \R^+$ introduced earlier, $b'$ has exactly one break 
point of increase $\bar a_0$ (with $\bar a_0 = -1/2$, in the identification 
$\torus \cong [-1/2, 1/2)$) and at most countably many break points 
of decrease $\{\bar a_i\}_{i \ge 1}$. In all other points of $\torus$, $b'$ is 
continuous. Also, $b'$ is decreasing (not necessarily strictly) on 
$\torus \setminus \{\bar a_0\}$. 

On the other hand, Lemma \ref{lem-xc} states that
$\map' = g_\theta \circ \arctan \circ \, b' \circ x_c$, where 
\begin{equation}
  g_\theta(\alpha) := \frac{\sin(\theta + \alpha)} {\sin(\theta - \alpha)}
\end{equation}
defines a strictly increasing continuous function $[-\alpham, \alpham] 
\into \R^+$ and $x_c : \torus \into \torus$ is the map defined in the 
statement of Lemma \ref{lem-xc}, which is easily seen to be an 
orientation-preserving homeomorphism. Therefore, there is a bijective 
correspondence between the break points of $b'$ and those of $\map'$,
such that the two functions have the same monotonicity properties 
between corresponding pairs of break points. This proves all the 
assertions of Lemma \ref{lem-non-reg}, except for the last one.

As for the last assertion, let us observe that, by definition, the function 
$\map'$ is not defined at its break points, so let us extend it to the 
whole of $\torus$ by setting $\map'(a_i) := \lim_{s \to a_i^+} 
\map'(s)$ (having employed the common abuse of notation whereby the 
extension has the same name as the extended function). It is evident 
that $\mathrm{Var}_{\torus \setminus \{a_0\}} (\log \map') = 
\Delta \log \map'(a_0)$, whence $\mathrm{Var}_{\torus} (\log \map') = 
2 \, \Delta \log \map'(a_0) < \infty$. (Here 
$\mathrm{Var}_{\torus \setminus \{a_0\}} ( \cdot )$ is the variation
of a real-valued function on the arc $\torus \setminus \{a_0\}$ and
$\mathrm{Var}_{\torus} ( \cdot )$ is the variation on the whole torus,
amounting to the former variation plus the variation at $a_0$.)
\qed

Now for the statements \emph{(ii)} of Theorem \ref{main-thm}. By 
Poincar\'e's dichotomy, if $\rho(\map) = p/q$ as in \emph{(ii)}, $\map$ 
has at least a periodic orbit of cardinality $q$, that is, $\map^q$ has at least 
$q$ fixed points. The previous proposition and the identity
\begin{equation} \label{main-30}
  (\map^q)'(x) = \prod_{k=0}^{q-1} \map'(\map^k(x))
\end{equation}
show that $(\map^q)'$ has at most $q$ break points of increase
(namely $\{ \map^{-k}(a_0) \}_{k=0}^{q-1}$, keeping in mind that some of 
these points may coincide), outside of which $\map^q$ is piecewise 
$C^1$ and concave. For the rest of this proof we refer to this property
as the `concavity property of $\map^q$'. Identifying $\torus$ with 
$[-1/2, 1/2)$, the graph of $\map^q$ can only look like one of the cases 
depicted in Fig.~\ref{fig3}. 

\newfig{fig3}{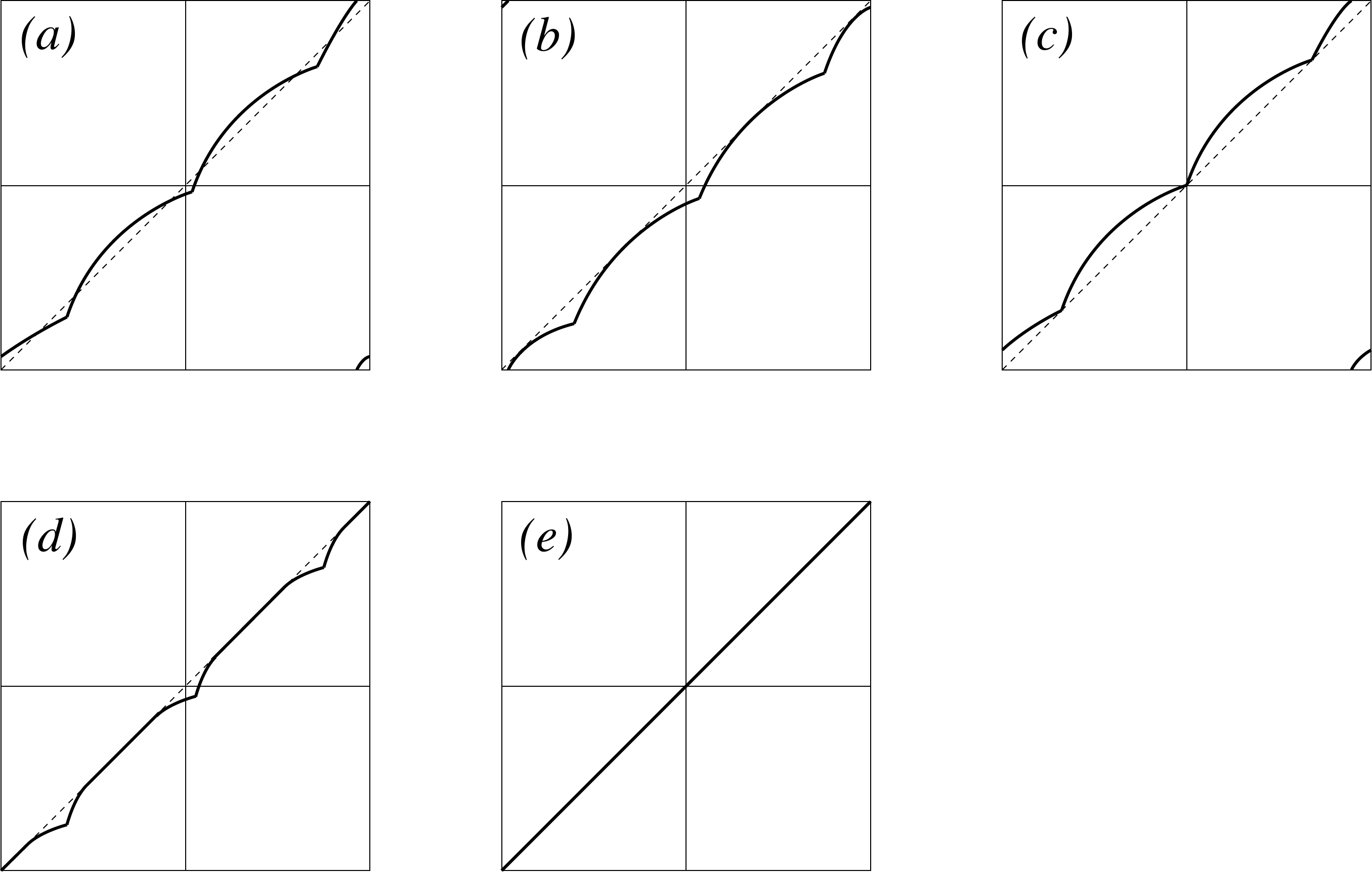}{12cm}{Possible graphs of the map 
$\map^q$, when a $q$-periodic orbit of $\map$ exists.}

In particular, the set of fixed points of $\map^q$ comprises at most $2q$, 
possibly degenerate, closed intervals $\perset_i^\pm$. By the concavity 
property of $\map^q$, the number of such intervals is $2q$ if, and only if, 
each $\perset_i^\pm$ is a singleton $\{ x_i^\pm \}$. In this case we label 
$x_i^-$, respectively, $x_i^+$, the repelling, respectively, attracting fixed 
points. Since, relative to the orientation of $\torus$, they alternate, the 
labelling can be done so that 
\begin{equation}
  x_0^- < x_0^+ < x_1^- < x_1^+ < \cdots < x_{q-1}^- < x_{q-1}^+ 
  < x_0^-. 
\end{equation}
Thus, by the second part of Poincar\'e's dichotomy, $\map^n(x_i^\pm) = 
x_{i+pm}^\pm$, for all $i \in \Z_q$ and $n \in \N$. This the case of
Fig.~\ref{fig3}\emph{(a)}. 

If the graph of $\map^q$ touches the bisectrix of the first and third
quadrants in a single point of abscissa $x$ then, if the graph touches 
``from below'' (Fig.~\ref{fig3}\emph{(b)}), $x$ is a fixed point of $\map^q$ 
which is repelling on the left and attracting on the right; if the graph touches 
``from above'' (Fig.~\ref{fig3}\emph{(c)}), $x$ is attracting on the left and 
repelling on the 
right. The same must therefore happen for all points of the $\map$-periodic 
orbit $\{ \map^i (x) \}_{i \in \Z_q}$, which must be distinct, by part
2 of Poincar\'e's dichotomy. The concavity property
of $\map^q$ shows that each of these points must belong to one and only
one concave part of $\map^q$, proving that there can be no other periodic
orbits of $\map$. Following their orientation on $\torus$, we label $x_0^-, 
x_1^-, \ldots, x_{q-1}^-$ the points of the aforementioned periodic orbit. 
In the case of left-repelling and right-attracting points 
(Fig.~\ref{fig3}\emph{(b)}), we also 
set $x_i^+ := x_i^-$; in the case of left-attracting and right-repelling points 
(Fig.~\ref{fig3}\emph{(c)}), we set $x_i^+ := x_{i+1}^-$. In either case, we 
finally denote $\perset_i^\pm := \{ x_i^\pm \}$. 

The above considerations prove all the assertions of Theorem 
\ref{main-thm}\emph{(ii)(a)}, when the sets $\perset_i^\pm$ are 
singletons. If $\perset^+ \cup \perset^-$ is not made up of isolated points, 
then it must include a closed interval, as we have seen. Take the largest 
interval (or, in case of a tie, one of the largest intervals) within 
$\perset^+ \cup \perset^-$. There can only be
two cases: either this closed interval is a proper subset of $\torus$,
in which case, without loss of generality, we denote it $\perset_0^+ := 
[x_0^L, x_0^R]$; or it is the whole of $\torus$, in which case we choose
any point $x_0^L \in \torus$ and set $x_0^R := \map(x_0^L)$, 
$\perset_0^+ := [x_0^L, x_0^R]$. In the former case, since $\map$ is a 
homeomorphism and preserves the property of being or not an 
$\map^q$-periodic point, we see that the sets 
$\{ \map^i (\perset_0^+) \}_{i \in \Z_q}$ are pairwise disjoint, closed 
intervals; they are also distinct, since all $f$-periodic orbits must have 
period $q$, cf.\ Fig.~\ref{fig3}\emph{(d)}. The concavity property of 
$\map^q$ shows that there can be no more periodic points. In the latter 
case, since the points $\{ \map^i (x_0^L) \}_{i \in \Z_q}$ are distinct, we 
have that the sets $[\map^i(x_0^L), \map^i(x_0^R)]$, $i \in \Z_q$, cover 
$\torus$ and intersect only at their endpoints. See Fig.~\ref{fig3}\emph{(e)}.

In either case, we denote $\perset_0^+, \perset_1^+, \ldots, 
\perset_{q-1}^+$ the intervals $\{ \map^i (\perset_0^+) \}_{i \in \Z_q}$, 
ordered according to the orientation of $\torus$, and 
$\perset_i^- := \perset_i^+$ for all $i$. The above facts prove all the 
claims of part \emph{(ii)(b)} of the theorem. In particular, 
Fig.~\ref{fig3}\emph{(d)} corresponds to case \emph{(ii)(b)(1)} and
Fig.~\ref{fig3}\emph{(e)} corresponds to case \emph{(ii)(b)(2)}.

Lastly, observe that if the function $b$ is strictly concave, (\ref{main-30}) 
shows that $\map^q$ is strictly concave in each of its concavity intervals, 
whence $\perset^+ \cup \perset^-$ can only contain isolated points, forcing 
the case \emph{(ii)(a)}.
\qed

\section{The trapezoid} 
\label{sect-trapez}

In this section we consider the case where $\btable$ is the rectangular 
trapezoid of Fig.~\ref{fig1}\emph{(b)}. We improve the general results of 
Theorem \ref{main-thm} and study the sets of parameters for which each 
of the only possible three cases happens: (1) there exist a global attractor 
and a global repellor; (2) all orbits are periodic with the same period; (3) all 
orbits are dense. 

We recall that the height of the trapezoid was fixed to 1/2 and that
in this case the two equal angles $\alpham = \alpha_m$ are denoted
$\alpha$. Let us fix the length of the shorter base to a certain
value $\ell>0$ and the direction of the flow to some angle $\theta \in 
(\alpha, \pi/2)$: this completely determines the homeomorphism 
$\map: \torus \into \torus$ introduced in Section \ref{sec-1ddyn}.
Lemma \ref{lem-xc} shows that the derivative $\map'$ only takes the 
values $\Lambda$ and $\Lambda^{-1}$, where
\begin{equation} \label{lamb}
  \Lambda := \frac{\sin(\theta+\alpha)}{\sin(\theta-\alpha)} >1.
\end{equation}
In order to simplify the ensuing computations, we impose the extra 
condition 
\begin{equation} \label{extra}
  \theta \ge \arctan(2\ell + \tan\alpha). 
\end{equation}  
The left out values of $\theta$ are qualitatively the same as the ones 
described below.
Recalling that the graph of $\map$ is symmetric around the bisectrix 
of the second and fourth quadrants of the square $[-1/2, 1/2)^2$,
we conclude that it must look like the one shown in Fig.~\ref{fig4}.
More in detail, using the notation of Lemma \ref{lem-non-reg}, let us 
call $a_0$ and $a_1$, respectively, the break points of increase and 
decrease of $\map'$. Simple calculations based on Fig.~\ref{fig2} 
show that 
\begin{equation} \label{a0-a1}
  a_0 = \frac{\tan \theta - 2\ell}{2 \tan \theta} \in \left( 0, \frac12 \right), 
  \qquad
  a_1 = -\frac{2\ell + \tan \alpha}{2 \tan \theta} \in \left( -\frac12, 0 \right),
\end{equation}
and $\map(a_j) = -a_j$, for $j=0,1$.

\newfig{fig4}{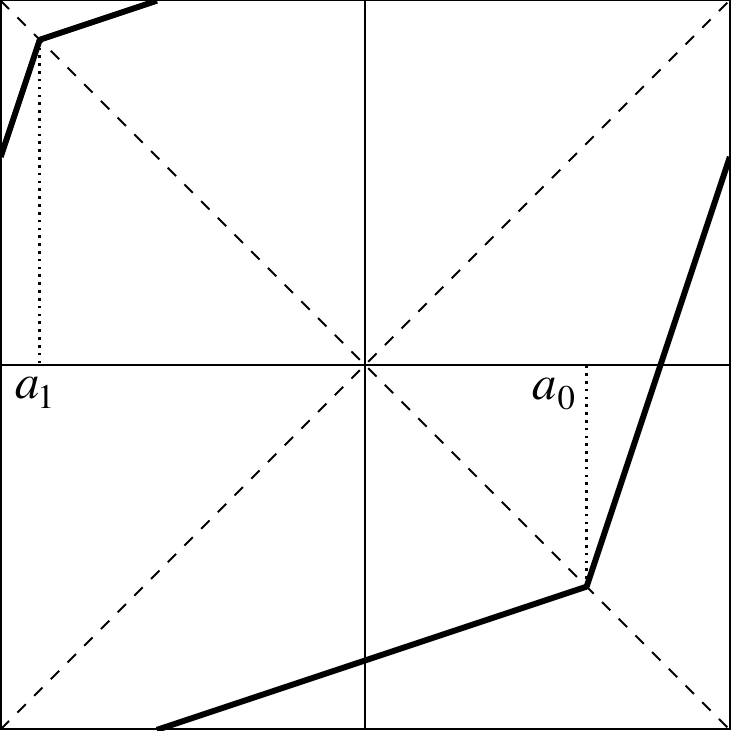}{5cm}{The graph of $\map$ in the case 
where $\btable$ is a rectangular trapezoid.}

In case our $\map$ has a rational rotation number, we can give more 
precise statements than Theorem \ref{main-thm}.

\begin{proposition} \label{prop-glob-reson}
  Let $f: \torus \into \torus$ be as defined above. Assume 
  $\rho(\map) = p/q$, with $p,q$ coprime positive integers, or $p=0$, 
  $q=1$. If $q$ is even then $\map^q=\id$, that is, all orbits are periodic 
  with primitive period $q$; moreover, $a_0$ and $a_1$ are in the same 
  periodic orbit. If $q$ is odd then, for any non-degenerate interval 
  $I \subset \torus$, $\map^q|_I \ne \id|_I$; also $a_0$ and $a_1$ are not 
  in the same periodic orbit.
\end{proposition}

\proof We first consider the case where $q$ is even. Let $x_0\in \torus$ 
be a periodic point with (primitive) period $q$ and denote by $x_0, x_1, 
\ldots, x_{q-1}$ the points of its orbit, labeled according to the 
orientation of $\torus$. 

By absurd, we assume that $a_0$ and $a_1$ are not periodic points. As 
explained in the proof of Theorem \ref{main-thm}, in this case $\map^q$ is 
piecewise linear with $2q$ break points of $(\map^q)'$, which correspond 
to the points in the backward orbits of $a_0$ and $a_1$ up to time $-q+1$. 
In other words, since $a_0$ and $a_1$ are not periodic points, the graph
of $\map^q$ falls in the case of Fig.~\ref{fig3}\emph{(a)}; more precisely, 
it is a 
polyline made up of $2q$ segments each of which crosses the bisectrix of 
the first and third quadrants. So, for all $i \in \Z_q$, there are exactly two 
break points between $x_i$ and $x_{i+1}$: one in the orbit of $a_1$ and 
one in the orbit of $a_0$ (recall that $i$ is a cyclic index, whence 
$x_q \equiv x_0$). 

For all points $x \in \torus$ such that $(\map^q)'(x)$ exists, we 
can write $(\map^q)'(x) = \Lambda^{n_+-n_-}$, where $n_+ = n_+(x)$ 
and $n_- = n_-(x)$ are the number of times the orbit of $x$ up to time 
$q-1$ falls in a set where $\map'$ is $\Lambda$ and $\Lambda^{-1}$, 
respectively. Clearly, $n_+ + n_-=q$, which is even, so $n_+-n_-$ is also 
even. As $x$ varies from $x_0$ to $x_1$, $(\map^q)'(x)$ varies at the 
two break points in the interval $(x_0,x_1)$: it is immediate to verify that 
passing through the point in the orbit of $a_0$ the effect is that $n_+$ 
increases by one and $n_-$ decreases by one, and the opposite 
happens passing through the point in the orbit of $a_1$. Therefore,
if $(\map^q)'(x_0) \ge \Lambda^2$, then $(\map^q)'(x) \ge 1$ for all 
$x \in (x_0,x_1]$ except for two break points. The last two inequalities
are in contradiction with $\map^q(x_0)=x_0$ and $\map^q(x_1)=x_1$.
Analogously, if $(\map^q)'(x_0) \le \Lambda^{-2}$, then
$(\map^q)'(x) \le 1$ for all $x \in [x_0,x_1]$ except for two break points, 
again a contradiction. Thus, $(\map^q)'(x_0) =1$. But this implies that 
$(\map^q)'(x) =1$ for all $x \in [x_0, \bar a)$, where $\bar a$ denotes the 
first break point to the right of $x_0$. All such points, then, including 
$\bar a$, are periodic, which is a contradiction, because $\bar a$ is in 
the orbit of $a_0$ or $a_1$.

We conclude that $a_0$ or $a_1$ must be periodic points, that is, 
we are in the cases \emph{(b), (c)} or \emph{(e)} of Fig.~\ref{fig3} (case 
\emph{(d)} is not achievable with a polyline of $2q$ segments). Let us 
assume for now that there exist periodic break points of $(\map^q)'$, 
that is, we are not in case \emph{(e)}. Let $x_0$ be a periodic break 
point. One can use essentially 
the same arguments as in the previous paragraph to prove that 
$(\map^q)'(x_0 + \eps) =1$ for all small $\eps>0$. (In fact, it is not 
possible that $x_0$ and $x_1$ are periodic points, $(\map^q)'(x) \ge 1$ 
for all but one $x $ in $(x_0, x_1)$, and 
$(\map^q)'(x_0+\eps) >1$ for a small $\eps>0$; the same goes for 
the analogous statement with reversed inequalities.) But this is 
incompatible with cases \emph{(b)} and \emph{(c)} of Fig.~\ref{fig3}. 

So only case \emph{(e)} is possible, i.e., $\map^q = \id$. This implies in 
particular that $a_0$ and $a_1$ belong to the same periodic orbit. In fact, 
if not, $(\map^q)'$  would have a positive jump at $a_0$, contradicting 
$(\map^q)' \equiv 1$.

Finally, let us consider the case where $q$ is odd. In this case 
$(\map^q)'(x) \ne 1$ for all $x$, since there is no choice for $n_+$ and 
$n_-$ defined above to satisfy $n_+=n_-$ and $n_++n_-=q$. One also
sees that $a_0$ and $a_1$ are not in the same periodic orbit, 
otherwise there would be no break points of $(\map^q)'$.
\qed

In the rest of the section it will often be convenient to study our
homeomorphisms by means of their lifts (cf.\ Section \ref{sec-1ddyn}). 
For a given $\map$ as described at the beginning of Section 
\ref{sect-trapez}, we consider the lift $F: \R \into \R$ uniquely defined by 
the conditions
\begin{align}
  & \map(x) = F(x) - \left\lfloor F(x) + \frac12 \right\rfloor, \quad \forall x \in
  \left[ -\frac12, \frac12 \right) ; \\
  & F(a_1) = -a_1.
  \label{ancoraggio}
\end{align}
This implies in particular that $F(a_0) = 1-a_0$ (recall that $\map(a_j) = 
-a_j$). The  last two equalities, together with the knowledge of the slope of 
$F$ before and after $a_0$ and $a_1$, cf.\ (\ref{lamb}), allow one to derive 
an expression for $F$:
\begin{equation} \label{lift-f}
  F(x) = \left\{ \begin{array}{ll} 
  \ds -a_1 + \Lambda^{-1} (x-a_1) , & \text{if } x \in [a_1, a_0] ; \\[4pt]
  1-a_0 + \Lambda (x-a_0), & \text{if } x \in [a_0, a_1+1] ; \\[4pt]
  F(x+k) - k,  & \text{if } x+k \in [a_1, a_1+1] \text{ with } k \in \Z.
\end{array} \right.
\end{equation}
An example of $F$ is shown in Fig.~\ref{fig5}.

\newfig{fig5}{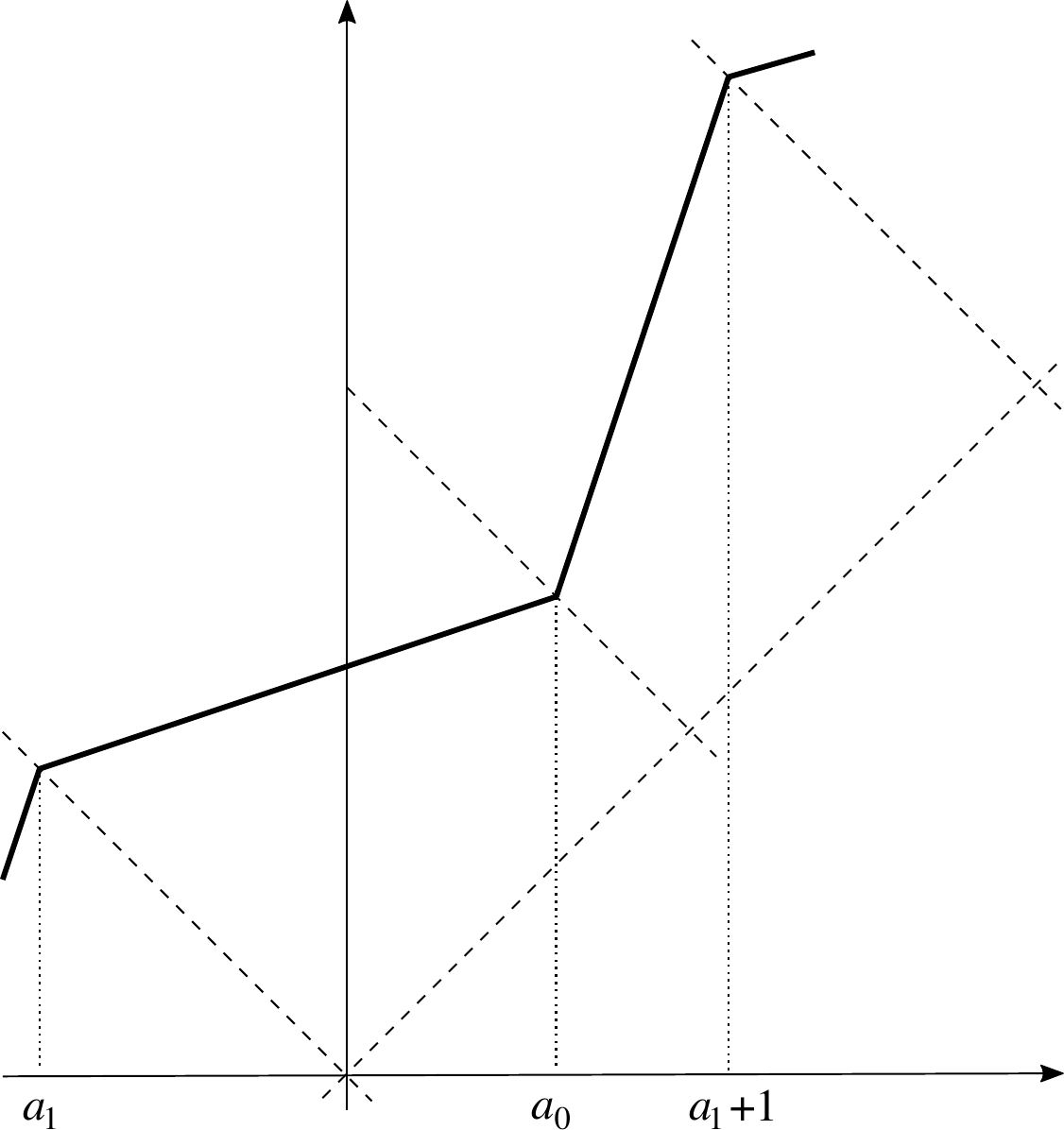}{8cm}{Graph of the lift $F$ corresponding to the 
map of Fig.~\ref{fig4}. The short dashed lines have equations $y=-x$, 
$y=-x+1$ and $y=-x+2$ in coordinates $(x,y)$.}

The main result of this part of the paper is that the typical $\map$
for a trapezoidal billiard has a rational rotation number. Here `typical' is 
intended in the strongest sense, that of measure theory. 

\begin{theorem} \label{thm-ae}
  Let $\map_{\ell, \alpha, \theta}: \torus \into \torus$ be the map introduced 
  earlier, for a given choice of the parameters $\ell>0$, $0 < \alpha < \theta < 
  \pi/2$. Then, for every $\ell$ and Lebesgue-a.e.\ $(\alpha, \theta)$,
  \begin{displaymath}
    \rho(\map_{\ell, \alpha, \theta}) \in \Q.
  \end{displaymath}
\end{theorem} 

The proof of this theorem, which uses more sophisticated techniques than 
the rest of the paper, is postponed to Appendix \ref{app-pwl2}. 

\medskip

A finer analysis of the properties of $\rho(\map)$, as $\map$ varies, 
comes from considering one-parameter families of maps. More specifically, 
we assume to be given a \emph{continuous strictly increasing family} 
$( \map_u )_{u \in [u_1, u_2]}$ of homeomorphisms of our type. This means 
that, for every $u \in  [u_1, u_2]$, there exists a lift $F_u$ of  $\map_u$ such 
that the family $(F_u)$ has the following properties:
\begin{itemize}
\item $u \mapsto F_u$ is continuous in the sup norm of $\R$;
\item for all $x \in \R$, $u \mapsto F_u(x)$ is strictly 
increasing;
\end{itemize}
cf.\ \cite[\S11.1]{kh}. A quick look at Fig.~\ref{fig5} convinces one that a
family of circle homeomorphisms in our class is continuous and strictly 
increasing if, and only if, the coordinates of the two break points $a_j$ of 
$\map'_u$ are continuous and strictly decreasing functions of $u$ (observe 
that the corner points of $F_u$ are constrained to lie on the lines $y = -x+k$, 
with $k \in \Z$). For example, employing again the notation 
$\map_{\ell, \alpha, \theta}$, it follows from (\ref{a0-a1}) that 
\begin{equation}
  ( \map_{\ell, \alpha, \theta} )_{\ell \in [\ell_1, \ell_2]}, \qquad
  ( \map_{\ell, \alpha, -\theta} )_{\theta \in [\theta_1, \theta_2]}, \qquad
  ( \map_{u, \kappa u, \theta} )_{u \in [u_1, u_2]},
\end{equation}
with $\kappa>0$, are continuous strictly increasing families whenever they
are defined.

\begin{theorem} \label{thm-tongues}
  Let $( \map_u )_{u \in [u_1, u_2]}$ be a continuous strictly increasing 
  family of circle homeomorphism of the type $\map_{\ell, \alpha, \theta}$, 
  with $\rho_1 := \rho(\map_{u_1}) < \rho_2 := \rho(\map_{u_2})$. Here 
  $\rho(\map_u)$ is intended as the actual r.h.s.\ of (\ref{rot-no}) with 
  $F = F_u$; in other words, $\rho(\map_u)$ can take all real values. Then 
  the function 
  \begin{displaymath}
    \rho: [u_1, u_2] \into \R, \quad \rho(u) := \rho(\map_u),
  \end{displaymath}
  is an increasing devil's staircase (namely, an increasing, continuous,  
  non-constant function that is constant on each interval of a family with 
  dense union). This implies in particular that each rotation number in 
  $[\rho_1, \rho_2]$ is realized by at least one $\map_u$. Moreover, 
  $\rho$ is constant on an interval if, and only if, the constant is 0 
  or $p/q$, with $p,q$ coprime integers and $q$ odd.
\end{theorem}

\proof It is a known fact that the function $\map \mapsto \rho(\map)$ is 
continuous and increasing in the space of orientation-preserving circle 
homeomorphisms. This means that the r.h.s.\ of (\ref{rot-no}), as an 
element of $\R$, is a continuous function of $F$, w.r.t the sup norm, which
 does not decrease if $F$ is increased pointwise; more details in 
 \cite[\S11.1]{kh}. Hence, $u \mapsto \rho(u)$ is continuous and increasing. 
 Moreover, by general results (see, e.g., \cite[Prop.~11.1.9]{kh}), it is strictly 
increasing at irrational values, i.e., if $\rho(u) \not \in\Q$ then, for all 
sufficiently close $u_- < u < u_+$, $\rho(u_-) < \rho(u) < \rho(u_+)$. 

The next proposition is a corollary of Proposition \ref{prop-glob-reson}.
We will use it to prove Theorem \ref{thm-tongues}, but the results are of 
independent interest as well.

\begin{proposition} \label{prop-odd-even}
  Denote by $Q_o$ the set of all rationals that can be written as $p/q$, 
  with $p$ and $q$ coprime and $q$ odd (by convention, $0=0/1 \in Q_o$). 
  Then: 
  \begin{itemize}
  \item[(i)] for all $r \in Q_o$, $\rho^{-1}(r)$ is a non-degenerate closed
    interval and, for all $u \in \rho^{-1}(r)$, $\map_u$ is not topologically 
    equivalent to a rotation;
      
  \item[(ii)] for all $r \in \R \setminus Q_o$, $\rho^{-1}(r)$ is a point 
    $u_r$ and $\map_{u_r}$ is topologically equivalent to the rotation by 
    $\rho(\map_{u_r})$.
  \end{itemize}
\end{proposition}

\proofof{Proposition \ref{prop-odd-even}} Let us first observe that, since
$F_u(x)$ is strictly increasing both in $u$ and $x$, the graphs of $F_u^q$
are strictly increasing in $u$, for all $q \in \Z^+$, i.e., 
$( \map_u^q )_{u \in [u_1, u_2]}$ is also a strictly increasing continuous 
family of circle homeomorphisms.

If, for a given $u_r$, $\rho(\map_{u_r}) = r \in Q_o$, Proposition 
\ref{prop-glob-reson} shows that $\map_{u_r}$ falls in the case 
\emph{(ii)(a)}  of Theorem \ref{main-thm} and so it cannot be conjugate 
to a rotation (an attractor exists). Moreover, by the continuity 
of $( \map_u^q )_u$, a left and/or right perturbation in $u$ preserves the
$q$ intersections between the graph of $\map_u^q$ and the bisectrix of 
the first and third quadrants in $[-1/2, 1/2)^2$. So, for such $u$, 
$\rho(\map_u) = r$. This ends the proof of part \emph{(i)}.

As for part \emph{(ii)}, we first consider $r \in \Q \setminus Q_o$ and
then $r \in \R \setminus \Q$. If $\rho(\map_u) = r \in \Q \setminus 
Q_o$, Proposition \ref{prop-glob-reson} gives $\map_u^q = \id$. By the 
strict monotonicity and $u$-continuity of $\map_u^q$, any arbitrarily small 
perturbation in $u$ 
will produce a positive distance between the graph of $\map_u^q$ and 
the bisectrix of the first and third quadrants, making $\rho(\map_u)$
arbitrarily close to an integer $p$, but unequal to $p$. Since it is a 
general fact that $\rho(\map^q) = q \rho(\map)$ mod 1, it follows that 
$\rho(\map_u)$ must vary, whence $\rho^{-1}(r) = \{u_r\}$. Moreover,
it is a general fact that if $\map^q = \id$ then $\map$ is topologically 
equivalent to the rotation by $p/q$, for some $p$. (Here is a sketch of its 
proof. By Poincar\'e duality, $\rho(\map) = p/q$, with $p$ coprime to $q$. 
Given a periodic orbit $\{ x_0, x_1, \cdots, x_{q-1} \}$, labeled according 
to the orientation of $\torus$, one defines $\phi |_{[x_0, x_1]}$ to be any 
orientation-preserving homeomorphism $[x_0, x_1] \into [0, 1/q]$. Then, 
for $x \in [x_i, x_{i+1}]$, $i \in \Z_q$, one chooses $n$ such that $np = i$ 
mod $q$ and defines $\phi(x) := \phi(\map^{-np}(x)) + i/q$. This gives a
homeomorphism $[x_i, x_{i+1}] \into [i/q, (i+1)/q]$. It is easy to verify
that, for $r \in [0,1)$, $\phi \circ \map \circ \phi^{-1}(r) = r + p/q$ mod 1.)

Finally, if $\rho(\map_{u_r}) = r \in \R \setminus \Q$, the fact that 
$\rho^{-1}(r) = \{u_r\}$ follows from the strict monotonicity of $\rho$ 
at irrational values, as recalled earlier, and the fact that $\map_{u_r}$ 
is topologically equivalent to the corresponding rotation follows from 
Theorem \ref{main-thm}\emph{(i)}.
\qed

A standard result of the theory of rotation numbers, cf.\ 
\cite[Prop.~11.1.11]{kh}, states that, if $( \map_u )_{u \in [u_1, u_2]}$ is a 
continuous monotonic family of orientation-preserving circle 
homeomorphisms such that $\rho(\map_u)$ is non-constant, and there 
exists a dense $S \subset \Q$ such that, whenever $\rho(\map_u) \in S$, 
$\map_u$ is not topologically equivalent to a rotation, then $u \mapsto 
\rho(\map_u)$ is a devil's staircase. We apply this result with $S := Q_o$ 
as in the statement of Proposition \ref{prop-odd-even}. This is obviously 
a dense subset of the rationals. Since $\rho$ is non-constant by
assumption, we conclude that $\rho$ is a devil's staircase, whose 
claimed properties have been proved earlier. This ends the proof of
Theorem \ref{thm-tongues}.
\qed

We give some final comments on the significance of the above results for 
the dynamics of internal-wave billiards in rectangular trapezoids.

The set of all rectangular trapezoids has 3 degrees of 
freedom, so the set of all internal-wave billiard flows with unit speed in a 
rectangular trapezoid has 4 degrees of freedom. On the other hand, as
for all billiard dynamics, a rescaling of the table gives rise to the same 
flow up to a rescaling of time. Moreover, an internal-wave billiard has 
another symmetry: a vertical or horizontal dilation of the table also leads 
to the same flow up to a rescaling of time. So the effective degrees of 
freedom in this problem are 2. We can choose to fix $\ell>0$ and represent 
the set of all flows as the parameter space $\mathcal{T} := \rset{(\alpha,\theta)} 
{0 < \alpha < \theta < \pi/2}$. 

Our results prove a few facts that were already observed, to a larger or 
smaller extent, in the physical literature:

\begin{enumerate}
\item \label{pt1} The parameter space is made up almost entirely by Arnol'd 
tongues. An Arnol'd tongue is the set of parameters for which the rotation 
number of the Poincar\'e map $\map$ assumes a given rational value, 
provided it has positive measure in parameter space. Theorem \ref{thm-ae} 
implies that almost every $(\alpha, \theta) \in \mathcal{T}$ is part of a 
tongue. Theorem \ref{thm-tongues}, applied to all families 
$( \map_{\ell, \alpha, -\theta} )_{\theta \in (-\pi/2, -\alpha)}$, for $\alpha \in
(0, \pi/2)$, shows that only rotation numbers with odd denominators (in 
simplest terms) possess a tongue. This was already observed in 
\cite{Maas1997}. In view of Proposition \ref{prop-glob-reson}, the above
result proves that, for almost all choices of the parameters, the billiard flow 
has a global attractor and a global repellor.

\item \label{pt2}The set of parameters $(\alpha, \theta)$ for which all 
trajectories are part of the same periodic beam is given by a countable 
number of smooth curves in $\mathcal{T}$. In fact, Proposition 
\ref{prop-glob-reson} shows that this case occurs if, and only if, 
$\rho(\map) = p/q$, with $q$ even (as a reduced fraction). It also shows
that such condition is equivalent to
\begin{equation} \label{eqs-case2}
  F_{\ell, \alpha, \theta}^q (a_0) = a_0+p, 
\end{equation}
where $F_{\ell, \alpha, \theta}$ is the lift of $\map = 
\map_{\ell, \alpha, \theta}$. By way of (\ref{a0-a1}) and (\ref{lift-f}), 
equation (\ref{eqs-case2}) is turned into a finite number of disjoint quadratic 
equations in the variables $\ell$, $\tan \alpha, \tan \theta$.

\item The set of parameters $(\alpha, \theta)$ for which all trajectories 
are dense is a zero-measure Cantor set, where a Cantor set is 
an uncountable set with no interior points. The nullity of the measure is a 
consequence of Theorem \ref{thm-ae}. The Cantor property comes from 
the proof of Theorem \ref{thm-tongues}, applied to the same families as 
in point \ref{pt1}: for all $\alpha$, no interval in the parameter $\theta$ may 
correspond only to irrational rotation numbers.

\item All rotation numbers in $[0,1)$ --- or, if we think of rotation numbers 
as taking values in the whole of $\R$, all rotation numbers in $\R^+$ --- 
are realized. In fact, it is easy to see that $F_{\ell, \alpha, \theta}$ 
converges to the identity, as $\theta \to \pi/2^-$, for all $\alpha$; whereas, 
for $\theta \to 0^+$, $F_{\ell, \theta/2, \theta} \approx T_\Delta$, where 
$T_\Delta: \R \into \R$ is the translation by $\Delta := 2\ell / \tan\theta$ (all 
convergences are in the sup norm of $\R$). With reference to point 
\ref{pt1} above, this proves in particular that \emph{all} rotation numbers 
with odd denominators (in simplest terms) possess Arnol'd tongues in 
$\mathcal{T}$.
\end{enumerate}

\begin{remark} \label{rem-mass-et-al}
The trapezoidal billiards considered in \cite{Maas1997} depend on 
two parameters, $d$ and $\tau$, where $\tau>0$ is the height of the 
trapezoid (in our case it is 1/2) and $d \in (-1,1)$ is such that the shorter 
base is $1+d$ (in our case it is $\ell$). In keeping with the fact that two 
parameters are enough to describe all internal-billiard flows, Maas 
\emph{et al} choose to fix the longer base of the trapezoid to be 2 and
$\theta = \pi/4$. The return map $\tilde \map$ they use is then defined on 
$\tilde \Sigma = [-\tau,\tau)$. It follows that $\tilde \map = 
\phi \circ \map \circ \phi^{-1}$, where $\phi : [-1/2, 1/2) \into [-\tau,\tau)$ is 
given by $\phi(x)=2\tau x$ and $\map$ is the map of Section 
\ref{sect-trapez} with parameters
\begin{equation} \label{chg-of-param}
  \ell=1+d, \qquad \alpha = \arctan(2 (1-d)) ,\qquad \theta = \arctan(2\tau).
\end{equation}
All results we prove for $\map$ hold for $\tilde \map$ with the 
corresponding parameters. 
In Fig.~\ref{fig6} we draw some Arnol'd 
tongues, computed rigorously with our methods, over the figure
\cite[Fig.~2]{Maas1997}, which describes simulations by Maas \emph{et al}.
More in detail, for some pairs $(p,q)$ of coprime integers, same-colored 
curves in Fig.~\ref{fig6} represent the boundary of the region in 
$(d,\tau)$-space where the rotation number equals $p/q$. The curves are
computed as follows. The proof of Theorem \ref{thm-tongues} shows that
the boundary of the mentioned region corresponds to the case where 
$a_0$ or $a_1$ are periodic orbits with rotation number $p/q$. These two 
occurrences are respectively equivalent to the two equations
\begin{equation} \label{eq-curve}
  F^q(a_j)=a_j+p, \qquad j=0,1. 
\end{equation}
Proposition \ref{prop-glob-reson} states that the two occurrences, and thus
the two equations, are distinct when $q$ is odd and the same when $q$
is even. As seen in point \ref{pt2} above, each equation amounts to a
finite number of disjoint equations in $\ell, \alpha, \theta$, which are
then turned into disjoint equations in $d,\tau$ via (\ref{chg-of-param}).
For example, in the case $(p,q) = (2,3)$, the lower and upper 
boundaries of the tongue are given respectively by
\begin{equation}
  \tau = \frac14 \! \left(d+5 + \sqrt{9d^2 + 10d + 17} \right) , \qquad
  \tau = \frac12 \! \left(d+2+\sqrt{d^2 + 8} \right).
\end{equation}
In the case $(p,q) = (1,4)$, to give an example with even denominator, the 
tongue degenerates to a single curve of equation
\begin{equation}
  \tau = d+3 + 2\sqrt{2+2d}.
\end{equation}
\end{remark}

\newfig{fig6}{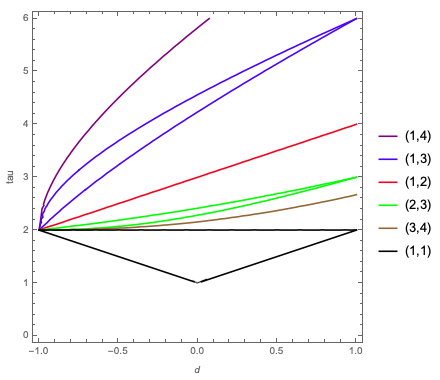}{10cm}{Arnol'd tongues in $(d,\tau)$-space. The 
figure shows the boundaries of some tongues relative to rotation numbers 
$p/q$, with $q$ odd, represented as curves of the same color. For $q$ even, 
each tongue degenerates to a single curve. The label of a tongue or a curve 
is the pair $(p,q)$. Compare with Fig.~2 of \cite{Maas1997}. }

\appendix

\section{Appendix: Piecewise linear circle homeomorphisms with two break 
points and dilation surfaces}
\label{app-pwl2}

In this appendix we prove Theorem \ref{thm-ae} by mapping our problem into
an analogous problem in the context of \emph{dilation surfaces} 
\cite{dfg, bfg, g}. Our arguments will be adaptations of arguments found in 
\cite{bfg}. For this reason, and also to avoid a very cumbersome appendix, we 
will not give every last detail of each proof, referring the reader to the 
explanations of \cite{bfg}.

\newfig{figa}{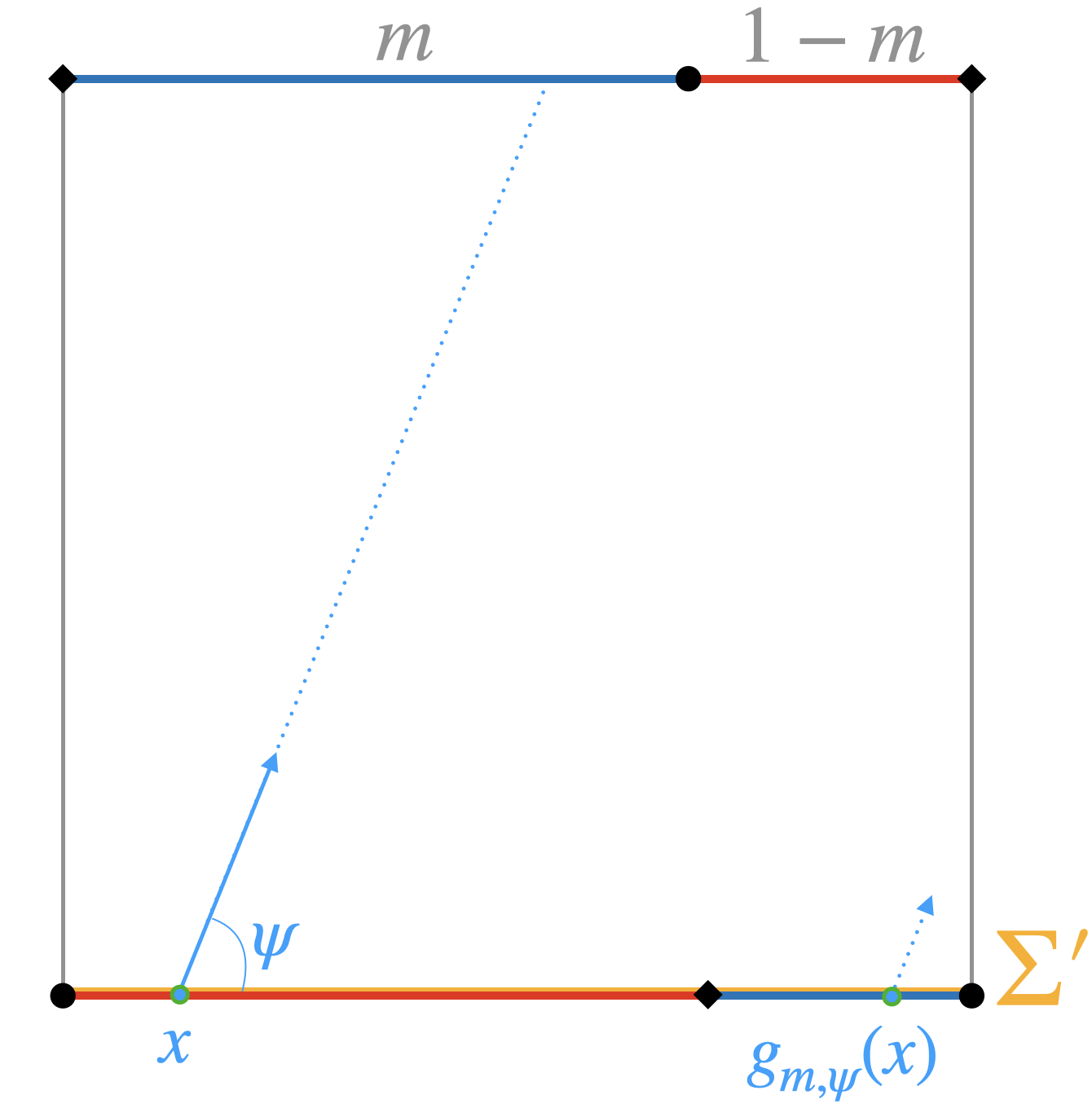}{5cm}{Dilation surface $\dilsurf \equiv \dilsurf_m$, 
  with $m \in [1/2, 1)$.  The gray vertical sides, of length 1, are identified by 
  means of a translation. The blue sides, of lengths $m$ and $1-m$ 
  respectively, are identified by means of an affine transformation. Same for 
  the red sides. Points marked with the same symbol are identified and
  represent a \emph{singular point} of $\dilsurf$. The curve $\Sigma'$, a 
  topological circle, is the bottom side of the square.}

The central idea is to view our maps $\map_{\ell, \alpha, \theta}$ as 
first-return maps for the linear flow on the dilation surface $\dilsurf$ 
illustrated in Fig.~\ref{figa}; see caption for a complete definition 
of $\dilsurf \equiv \dilsurf_m$, for $m \in [1/2, 1)$. More in detail, in view of 
Fig.~\ref{figa}, consider the flow on $\dilsurf$ defined by the constant vector 
field of direction $\psi \in (0, \pi/2]$. We refer to it as the \emph{linear flow} 
of direction $\psi$. The curve $\Sigma'$ has 
the topology of a circle and is clearly a cross-section for the linear flow.
We call $g_{m, \psi}: \Sigma' \into \Sigma'$ the corresponding Poincar\'e
map, where $m \in [1/2, 1)$ is the parameter indicated in Fig.~\ref{figa}. 

Now, to the triple of parameters $(\ell, \alpha, \theta)$, with $\ell>0$ and 
$0 < \alpha < \theta \le \pi/2$, we associate the pair $(m, \psi) \in
[1/2, 1) \times (0, \pi/2]$ with
\begin{equation} \label{m-psi}
  m := \frac12 + \frac{\tan \alpha}{2\tan \theta} , 
  \qquad
  \psi := \arctan \left[ \left( \frac{2\ell} {\tan \theta} +m-1 \right)^{-1} 
  \right] .
\end{equation}
We claim that, up to conjugation, $g_{m, \psi}$ acts on $\Sigma'$ as 
$\map_{\ell, \alpha, \theta}$ acts on $\Sigma$. First observe that, by 
(\ref{m-psi}) and (\ref{a0-a1}), $m = a_0-a_1$. This is actually the reason 
behind the definition of $m$. Also, cf.\ (\ref{lamb}) and Fig.~\ref{fig4},
\begin{equation} 
  \Lambda = \frac{a_0-a_1}{1-a_0+a_1} = \frac m {1-m} .
\end{equation}
Redrawing Fig.~\ref{fig2} for the case where $\ubtable$ is a hexagon 
(corresponding to $\btable$ being our rectangular trapezoid) shows 
that, up to a rotation of $\Sigma$, the action of $\map_{\ell, \alpha, \theta}$ is 
completely determined by the fact that an arc of length $1-m$ (the arc 
$[a_0, a_1]$ in Fig.~\ref{fig4}) is uniformly expanded by a factor $\Lambda = 
m/(1-m)$ and its left endpoint is moved to the right by quantity 
$2\ell / \tan \theta$ (here left and right are relative to the orientation of 
$\Sigma$). The complementary arc, of length $m$, follows suit: it is 
uniformly contracted by a factor $\Lambda^{-1}$ and its left endpoint is 
moved to the right by a quantity $(2\ell + \tan \alpha) / \tan \theta$.
On the other hand, by construction, $g_{m, \psi}$ uniformly expands an
arc of length $1-m$ by a factor $m/(1-m) = \Lambda$ and moves its left 
endpoint to the right by a quantity $1/\tan \psi + 1-m$. Again, the 
complementary arc, of length $m$, follows suit: it is uniformly contracted by a 
factor $\Lambda^{-1}$ and its left endpoint is moved to the right by a quantity 
$1/\tan \psi + m$. So, for any isometric identification $\Sigma \cong \Sigma'$,
$\map_{\ell, \alpha, \theta}$ coincides with $g_{m, \psi}$ up to a conjugation
by rotation if, and only if,
\begin{equation} 
  \frac{2\ell}{\tan \theta} = \frac1 {\tan \psi} + 1-m,
\end{equation}
which is exactly the second relation of (\ref{m-psi}).

In the remainder, we will simplify the notation by introducing the new 
parameter
\begin{equation} \label{def-s}
  s := \frac1 {\tan \psi} = \frac{2\ell} {\tan \theta} + m-1 = \frac{4\ell + \tan \alpha} 
  {2\tan \theta} - \frac12
\end{equation} 
and rewriting $g_s := g_{m, \psi}: \torus \into \torus$ (via the natural 
identification $\Sigma' \cong \torus$). The latter change of notation is 
convenient because we are interested in fixing $\ell, m$ and varying $s$. The 
main goal of this appendix is to show that the rotation number of $g_s$ is 
irrational for a.e.\ $s$.

\begin{theorem} \label{thm-app}
  For all $m \in [1/2, 1)$ and Lebesgue-a.e.\ $s \in [0,+\infty)$, $\rho(g_s) 
  \in \Q$.
\end{theorem}

\proof Here is where non-trivial geometric tools for dilation surfaces come into 
play. We will see that it is useful to consider the action that the group of 
orientation-preserving affine diffeomorphisms of $\dilsurf$ has on the directional 
foliations of that surface. A \emph{directional foliation} of $\dilsurf$ is simply the 
collection of all the trajectories of a given linear flow on $\dilsurf$. For such an 
action, obviously, the translational part of an affine diffeomorphism is irrelevant, as 
is any scaling factor in the linear part. In other words, the relevant information is 
stored in $V_{\dilsurf}$, the \emph{Veech group} of $\dilsurf$, which is defined to 
be the group of the linear parts of all orientation-preserving affine diffeomorphisms 
of $\dilsurf$, modulo multiplicative factors. Having chosen a reference system for 
$\dilsurf$ (in the present case, the one implicit in Fig.~\ref{figa}), $V_{\dilsurf}$ 
is regarded as a subgroup of $SL(2,\R)$.

It turns out, cf.~\cite[Thm.~4]{bfg}, that $V_{\dilsurf}$ is the group generated by 
the matrices
\begin{equation}
  A = \begin{pmatrix} 1 & 1 \\ 0 & 1 \end{pmatrix} , \qquad 
  B = \begin{pmatrix} 1 & 0 \\ 1/(m - m^2) & 1 \end{pmatrix} , \qquad 
  -I = \begin{pmatrix} -1 & 0 \\ 0 & -1\end{pmatrix}.
\end{equation} 
In the interest of an accessible exposition we show, by means of Figs~\ref{fige} 
and \ref{figf}, that $A,B \in V_{\dilsurf}$ (the fact that $-I \in V_{\dilsurf}$ is 
obvious). We do not prove, however, that the group generated by $A, B, -I$ is the 
entire Veech group, because the forthcoming arguments work as well if one 
replaces $V_{\dilsurf}$ with $\langle A, B, -I \rangle$.

\newfig{fige}{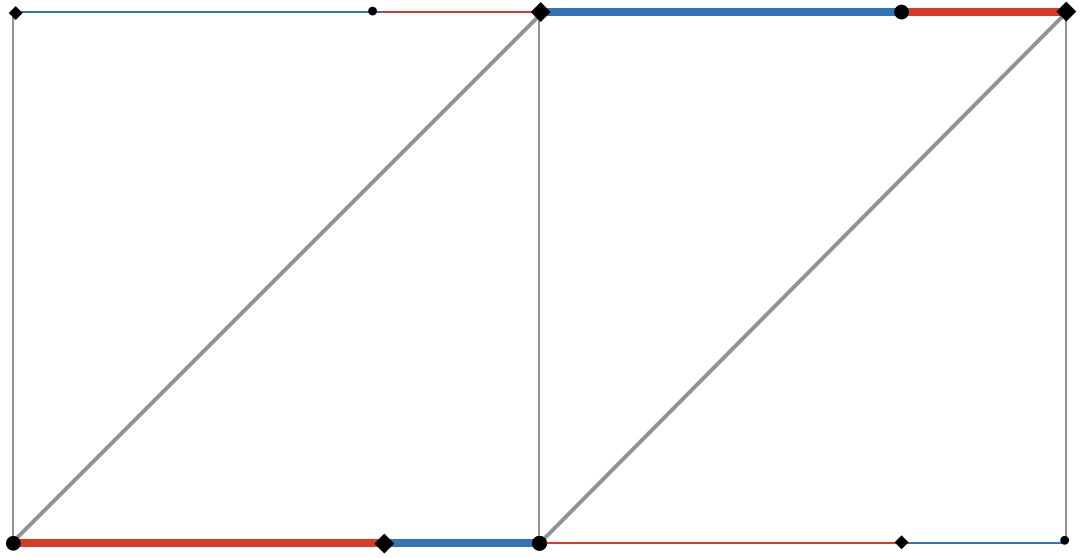}{7.4cm}{The action of $A$ on $\dilsurf$. More precisely, by
  way of a two-fold cover of $\dilsurf$, the figure shows the action of 
  an orientation-preserving affine diffeomorphism of $\dilsurf$ 
  whose linear part is $A$.}

\newfig{figf}{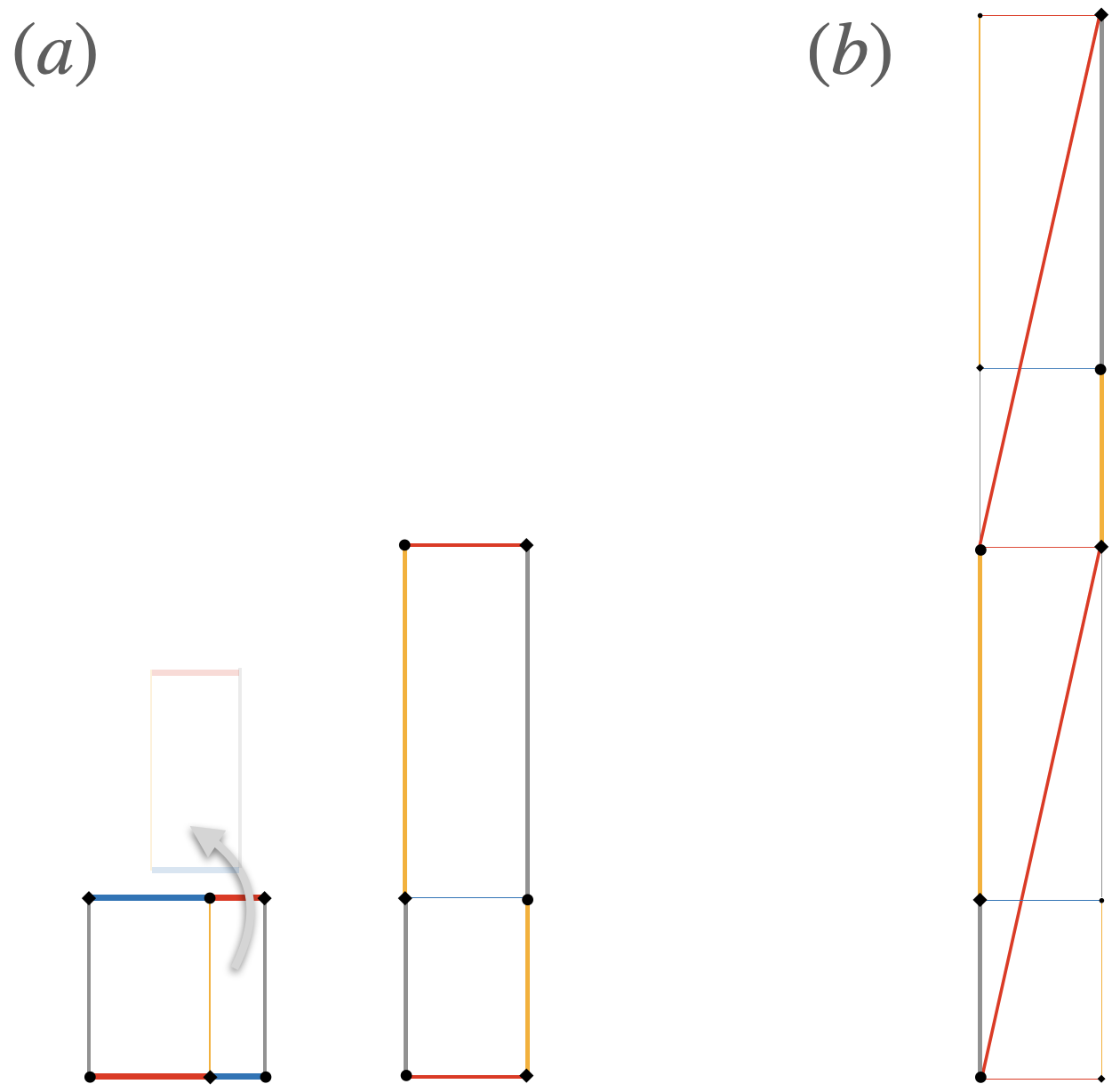}{10cm}{The action of $B$ on $\dilsurf$. As in Fig.~\ref{fige}, 
  part \emph{(b)} shows the action of an orientation-preserving 
  affine diffeomorphism of $\dilsurf$ whose linear part is $B$, but here we are 
  using another planar representation of $\dilsurf$, constructed in part
  \emph{(a)}.}

We parametrize the directional foliations by means of the variable $s = 1/\tan \psi$,
where we extend the domain of $\psi$ to $[0, \pi)$. This means that we regard 
$s$ as a projective variable in $\R P^1 \cong \R \cup \{\infty\}$. Notice that, in
so doing, we are only considering directional foliations with a non-negative 
vertical component or, more truthfully, we are considering foliations up to their 
orientation. To see how an element
\begin{equation}
  M = \begin{pmatrix} a & b \\ c & d \end{pmatrix} 
\end{equation} 
of $V_{\dilsurf}$ acts on the foliation parametrized by $s$, consider a vector 
$v = (r, r\tan\psi) = (r, r/s)$, $r \in \R \setminus \{0\}$, in the direction of such 
foliation (modulo orientation). In the coordinate $s$, the action of $M$ is clearly
\begin{equation} \label{action-on-s}
  s = \frac r {r/s} \mapsto \frac{ar + b(r/s)} {cr + d(r/s)} = \frac{as + b}{cs + d}.
\end{equation} 
With a slight abuse of notation, we denote the above r.h.s.\ by $M(s) \in \R P^1$. 
Obviously, $s \mapsto M(s)$ is the M\"obius transformation given by the matrix 
$M$, regarded as an element of $PSL(2,\R)$ (with the standard abuse). 

We call $\Gamma_{\dilsurf}$ the projection of $V_{\dilsurf}$ into $PSL(2,\R)$. 
Thus, $\Gamma_{\dilsurf}$ is the Fuchsian group generated by $A$ and $B$.
In this spirit, for $M \in \Gamma_{\dilsurf}$, we extend the definition of $M(z) := 
(az+b)/(cz+d)$ to $z \in \HHH$, the hyperbolic upper half plane. In this way, 
$\Gamma_{\dilsurf}$ acts on $\HHH$, so its action on $\R P^1$ can be studied 
as its action on the ``boundary at infinity'' of $\HHH$.

It is shown in \cite[\S3.3]{bfg} that the maximal set of $\R P^1$ where the action
of $\Gamma_{\dilsurf}$ is properly discontinuous is a set $D_\Gamma$ of full
Lebesgue measure and that $D_\Gamma / \Gamma_{\dilsurf}$ has the topology
of a circle. A corresponding fundamental domain $J_\Gamma$ can be found by 
taking the ``boundary at infinity'' of a fundamental domain $\mathcal{F}_\Gamma$ 
for the action of $\Gamma_{\dilsurf}$ on $\HHH$. An example of 
$\mathcal{F}_\Gamma$, calculated by means of $A$ and $B$, is depicted in 
Fig.~\ref{figg}, together with its corresponding 
\begin{equation}
  J_\Gamma = \left[ \frac1 {-2 + 1/(m-m^2)} \, ,\, \frac12 \right] .
\end{equation} 

\newfig{figg}{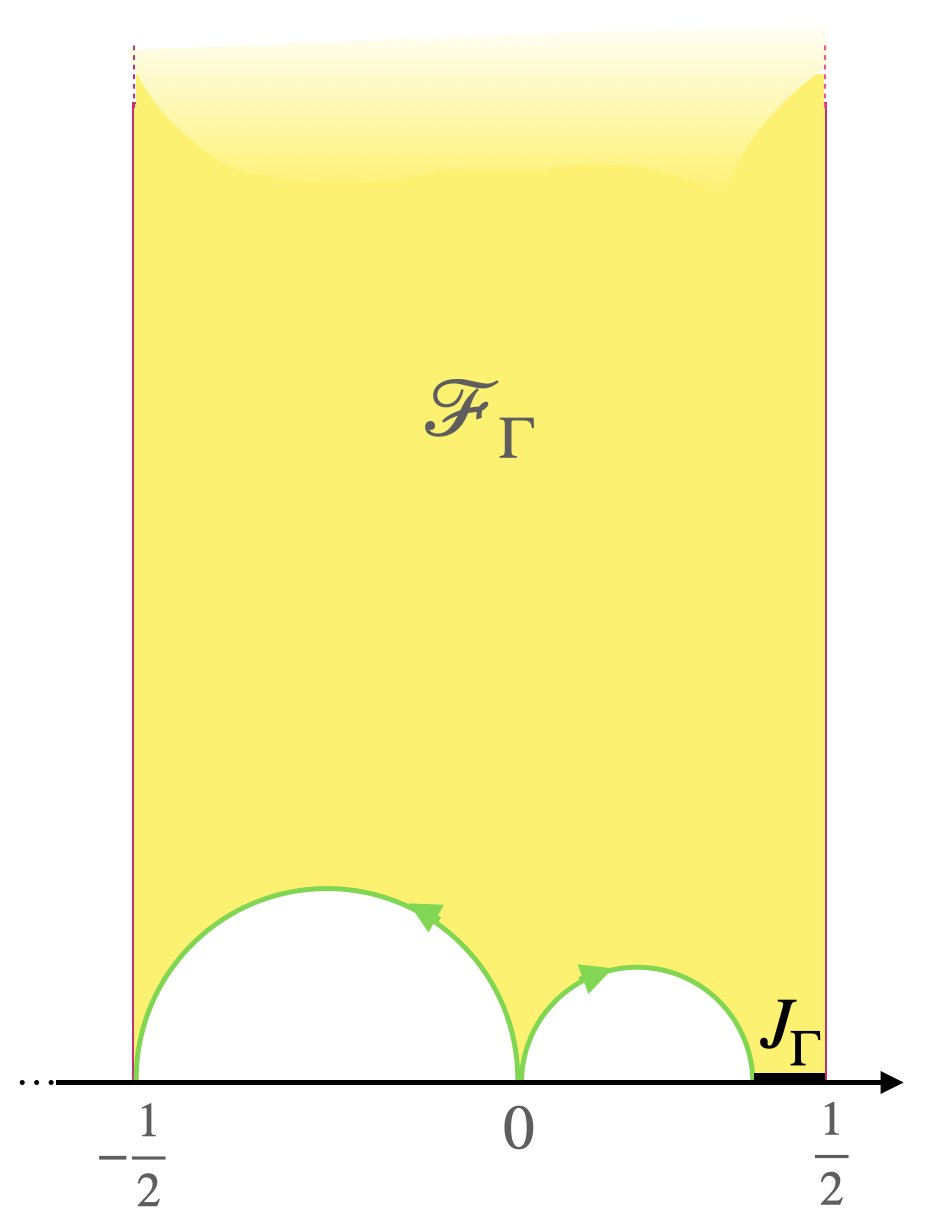}{6cm}{The fundamental domains $\mathcal{F}_\Gamma$ 
  and $J_\Gamma$.}

\begin{lemma} \label{lemma-app1}
  Every map $g_s$ with $s \in J_\Gamma$ has a fixed point.
\end{lemma} 

\proofof{Lemma \ref{lemma-app1}} The statement is true for a larger interval of 
directions, namely, $[1-m, m] \supset J_\Gamma$. In fact, as in Fig.~\ref{figl},
one can draw the cone of half-lines that homothetically projects that upper red 
segment into the lower red segment of $\dilsurf$. The directions defined by this
cone, in the projective coordinate $s$, are exactly $[1-m, m]$. It is easily 
seen that the intersection of every half-line of the cone with $\dilsurf$ is a closed
trajectory for the linear flow defined by the corresponding direction $s$, with the
property that the trajectory closes up at the first return to $\Sigma'$. In other words,
the trajectory corresponds to a fixed point of $g_s$.
\qed

\newfig{figl}{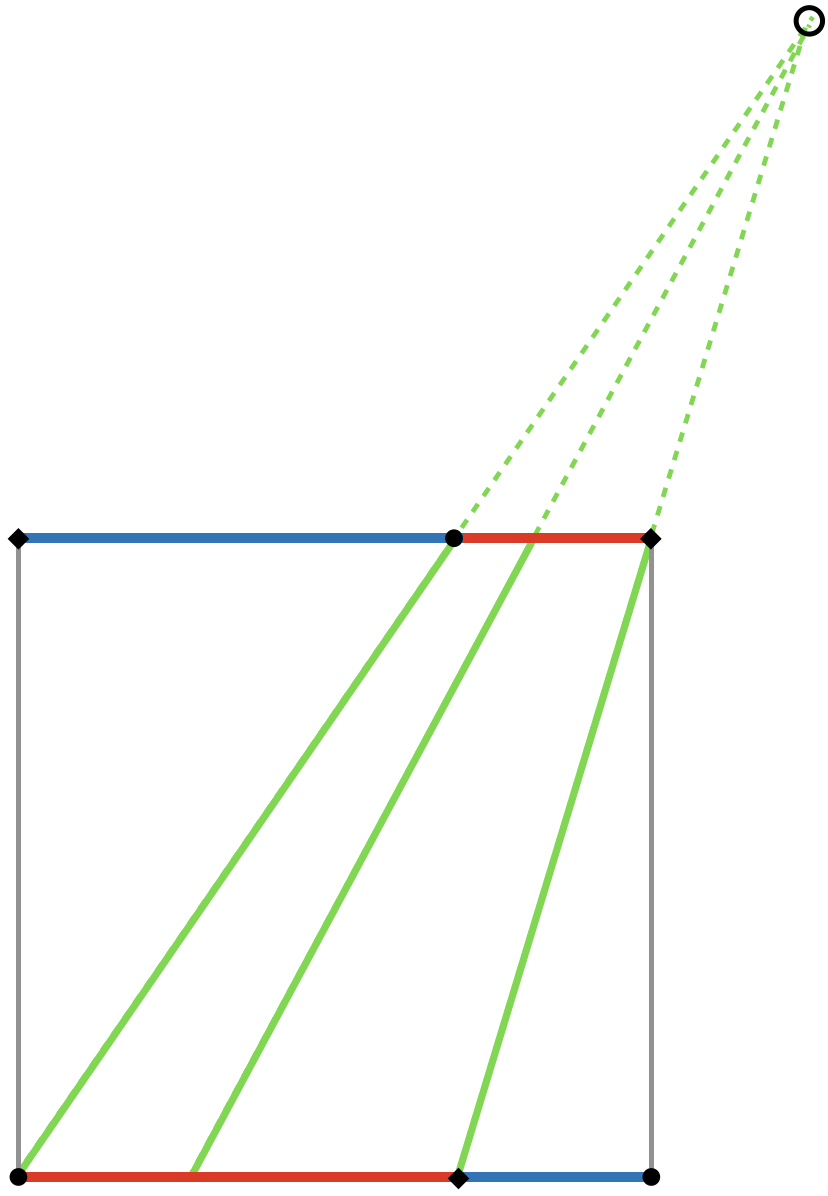}{5cm}{Closed leaves for $s \in [1-m, m]$.}

\begin{lemma} \label{lemma-app2}
  For $s \in \R P^1$ and $M \in \Gamma_{\dilsurf}$, $g_s$ has a periodic orbit
  if, and only if, $g_{M(s)}$ has a periodic orbit.
\end{lemma} 

\proofof{Lemma \ref{lemma-app2}} The proof is immediate: $g_s$ has a periodic 
orbit if, and only, if the foliation labeled by $s$ has a closed leaf. Now, 
$M \in \Gamma_{\dilsurf}$ means that there exists an orientation-preserving 
affine diffeomorphism $\tau_M$ of $\dilsurf$ whose linear part, up to factors, is 
$M$. By definition of $M(s)$, cf.\ (\ref{action-on-s}), $\tau_M$ transform the
foliation labeled by $s$ into the foliation labeled by $M(s)$ and the closed leaf
for the former into a closed leaf for the latter, and viceversa.
\qed

\smallskip

We are now in a position to quickly finish the proof of Theorem \ref{thm-app}.
By Poincar\'e's classical theory (see proof of Theorem \ref{main-thm}), 
$\rho(g_s) \in \Q$ if, and only, if $g_s$ has a periodic orbit. The previous 
statements imply that, for every $s \in D_\Gamma$, there exists an $M \in 
\Gamma_{\dilsurf}$ such that $M(s) \in J_\Gamma$. By Lemmas \ref{lemma-app1}
and \ref{lemma-app2}, $\rho(g_s) \in \Q$. Restricting $s$ to $D_\Gamma \cap
[0, +\infty)$ yields the theorem.
\qed

\proofof{Theorem \ref{thm-ae}} Take $\ell > 0$ constant. Fixing $m \in [1/2, 1)$ 
means fixing the ratio $\tan \alpha / \tan \theta = 2m-1 \in [0,1)$, cf.\ the first 
relation of (\ref{m-psi}). By (\ref{def-s}), an almost sure property in  terms of $s$ 
translates in an almost sure property in terms of $\tan \theta$. So, in every 
half-line of the strip 
\begin{equation}
  \left( \tan \theta, \frac{\tan \alpha} {\tan \theta} \right) \in (0,+\infty) \times [0,1),
\end{equation}
the property $\rho(\map_{\ell, \alpha, \theta}) \in \Q$ occurs Lebesgue-almost 
everywhere. Since the transformation $(\alpha, \theta) \mapsto 
(\tan\theta, \tan \alpha / \tan \theta)$ has positive Jacobian, the same property 
holds Lebesgue-almost everywhere in the triangle $\mathcal{T} := 
\rset{(\alpha,\theta)} {0 < \alpha < \theta < \pi/2}$.
\qed

\begin{remark}
  Thinking of the directional foliations of $\dilsurf$ as oriented upwards, observe 
  that the proof of Lemma \ref{lemma-app1} provides a \emph{repelling} fixed point 
  of $g_s$. One can improve the lemma by also drawing the cone joining the two 
  blue sides. 
  The directions defined by this cone are again $[1-m, m]$, but all half-lines in it 
  produce an attractive closed trajectory. Therefore, every map $g_s$ with $s \in 
  [1-m, m]$ has an attractive and a repelling fixed point.
  This allows us to improve the statement of Theorem \ref{thm-ae} as well. In fact,
  the proof of Theorem \ref{thm-app} now shows that, for a.e.\ $s \in \R P^1$, 
  $g_s$ has an attracting and a repelling periodic orbit. Thus, by way
  of Proposition \ref{prop-glob-reson}, the proof of Theorem \ref{thm-ae} implies 
  that, for all $\ell > 0$ and a.e.\ $(\alpha, \theta)$, $\rho(\map_{\ell, \alpha, \theta})$ 
  can be written in the form $p/q$, with $q$ odd. This result is also
  established by other methods at the end of Section \ref{sect-trapez}.
\end{remark}

\end{document}